\magnification=\magstep1
\baselineskip=18pt
\parskip=6pt
\input amstex
\documentstyle{amsppt}
\vsize=9truein
\loadmsbm
\loadeufm
\NoBlackBoxes
\topmatter
\title Local Factors, Reciprocity and Vinberg Monoids\endtitle
\author F.~Shahidi\footnote"*"{Partially supported by the NSF grant.\hfil\break} \endauthor

\loadmsbm
\loadeufm
\UseAMSsymbols

\def\ad{\text{ad}}
\def\Gal{\text{Gal}}
\def\der{\text{der}}
\def\alg{\text{alg}}
\def\sc{\text{sc}}
\def\on{\overline n}

\def\Sym{\text{Sym}}
\def\Int{\text{Int}}
\def\stand{\text{stand}}
\def\char{\text{char}}
\def\stan{\text{stan}}
\def\trace{\text{trace}}
\def\std{\text{std}}
\def\End{\text{End}}
\def\ff{\frak p}
\def\bC{\Bbb C}
\def\bQ{\Bbb Q}
\def\bA{\Bbb A}
\def\bN{\Bbb N}
\def\bZ{\Bbb Z}
\def\bS{\Bbb S}
\def\bG{\Bbb G}
\def\Aut{\text{Aut}}
\def\Card{\text{Card}}
\def\Ind{\text{Ind}}
\def\oK{\overline K}
\def\ok{\overline k}
\def\var{\varepsilon}
\def\opsi{\overline\psi}
\def\oN{\overline N}
\def\sH{\Cal H}
\def\sM{\Cal M}
\def\sS{\Cal S}

\def\check{\vee}

\centerline{\bf Local Factors, Reciprocity and Vinberg Monoids}
\bigskip
\centerline{F.~Shahidi\footnote"*"{Partially supported by the NSF grant DMS--1500759}}
\bigskip
\centerline{\bf Introduction}

This article addresses the problem of existence of local factors, i.e., the root numbers and $L$--functions attached to representations of reductive groups over local fields and irreducible finite dimensional representations of their $L$--groups, as well as their equality with those of Artin factors through the local Langlands correspondence (LLC).

There are several methods for defining the local factors, notably the Rankin--Selberg [JPSS,S1,S2], Langlands--Shahidi [Sh5,Sh8], and now the approach of Braverman--Kazhdan [BK] which generalizes that of Tate [T1] and Godement--Jacquet [GJ] in the case of the standard $L$--functions for $GL_n$ to a fairly general setting, and has attracted the attention of those interested in Beyond Endoscopy [L5] through the conjectural presence of a Fourier transform and its Poisson summation formula [FLN].

On the other hand the recent work of Laurent Lafforgue [La] explains how a global theory of $L$--functions can produce a Fourier transform and Poisson summation formula, presenting the global functional equation as a Poisson summation, thus, in some sense putting different theory of $L$--functions in the same footing, and connecting them to the theory of basic functions and functoriality [BNS,ChN].
See Remark 9.4.16.

The approach in [BK] is quite general and as explained fairly nicely by Ngo [N1,N2], what replaces the vector space $M_n(F)$ of the standard $L$--function of Godement--Jacquet [GJ] when the standard representation is changed to an arbitrary irreducible representation of $^LG$ with highest weight $\lambda$, is a monoid $M^\lambda$ attached through the Vinberg's theory of universal monoids [V].
The group $G^\lambda$ is then the group of units of $M^\lambda$ (cf.~Section 9).

In this paper we address two problems.
The first is to establish the equality of local (analytic) factors defined by a theory of $L$--functions attached to an arbitrary irreducible finite dimensional representation $r$ of $GL_n(\bC)$, the $L$--group of $GL_n$, satisfying certain natural axioms, with those of Artin factors through (LLC), established for $GL_n$ in [HT,He1,Sch].
The axioms are precisely those meant and needed for local--global compatibility for LLC.

They are in fact natural generalization of the work in [CST], joint with J.~Cogdell and T.-L.~Tsai, where the equality of the factors for $r=\Lambda^2$ or $\Sym^2$, the exterior and symmetric square representations of $GL_n(\bC)$ are proved.
As it is shown here the arguments of [CST] are robust enough to be extended to an arbitrary $r$ in the case of $GL_n$.

The work in [CST] is reviewed in Sections 3--6, after which we address the general $r$ in Section 8.
The new ingredient is the generalization of the axiom of ``multiplicativity'' which requires the use of Schur functors and Young symmetrizers [FH] through the Littlewood--Richardson rule which we review in (8.3.1)--(8.3.32).
We then prove the equality of arithmetic (Artin) and analytic factors through LLC using this generalized version of multiplicativity ``$M$'' stated in (8.3.34).
The final theorem is stated as Theorem 8.3.50.

The representations $r$ is given by a partition $\nu$ of length $|\nu|$ which defines the Schur functor $\bS_\nu$ which acts on every underlying space $\bC^n$, defining a representation of $GL_n(\bC)$ for each $n$.
For a general representation of the Weil group, following Harris [H,He2], the proof of (8.3.50) requires the validity of axioms for all $n$.
But when the representation is monomial and in particular for all tame cases, the theorem does not require any other group which shows how much more difficult the wild cases are.

We should point out that even if the representation of the Weil group is monomial, its composition with an arbitrary $r$ may not necessarily decompose to a sum of monomial representations and thus one will need to deal with non--monomial representations at any rate for a general $r$.
This is an important and restrictive observation and shows the difference of the $\var$--factor for a general $r$ versus the cases we have so far confronted such as $r=\Lambda^2$ and $r=\Sym^2$.
Our local--global compatibility axioms and approach allows us to virtually circumvent this difference and difficulty.

Moreover, in view of the recent progress (e.g.~[JLi]) in establishing Jacquet's conjecture on local converse theorems for $GL(n)$, the importance of stability under highly ramified twists becomes clearly evident whenever one attempts to apply converse theorems to establish cases of functoriality, as the global version comes within reach.
We expect that our approach and progress towards stability will play a crucial role in establishing such results.

The last part of this paper, Section 9, is devoted to a computation of the group $G^\lambda$, invertible elements of the monoid $M^\lambda$ following the exposition by Ngo [N1,N2] of Vinberg's theory of universal monoids, in general, and its connection to Langlands--Shahidi method and beyond.
More precisely, we note (Proposition 9.3.12 --- communicated to us by Steve Miller) that for every irreducible finite dimensional complex representation $\rho$ of a simply connected complex Lie group $G$ for which $G/G_{\der}$ is one dimensional, there exists a complex Kac--Moody group $H$ and a maximal parabolic subgroup $P\subset H$ with a Levi decomposition $P=LN$ with $L_{\der}=G_{\der}$ such that $\rho$ appears in the adjoint action of $L$ on Lie$(N)$.
This includes all the cases appearing in Langlands--Shahidi method in finite dimensional cases of $H$.

What we observe in this section as our Proposition 9.3.11, is that we can in fact choose the pair $(H,L)$ such that $L\simeq\hat G^\lambda$, the $L$--group of $G^\lambda$ with $G^\lambda$ the group of units of $M^\lambda$ as above.

As examples we look at $GL_n$ and the cases $\rho=\Sym^m$ and $\rho=\Lambda^m$ for which $\lambda=m\delta_1$ and $\lambda=\delta_m$, respectively, and compute $G^\lambda$ in all cases when $m$ is prime as well as the cases of Rankin products for $GL_n\times GL_m$.
They all agree with what happens in the lists in [L1,Sh4,Sh8].
Finally, in the case of $GL_2$ we also compute the monoids $M^\lambda$ attached to all symmetric powers of $GL_2(\bC)$ in paragraph (9.3.18).

We conclude Section 9 by showing that the local coefficients [Sh1] which allow us to define the root numbers and $L$--functions from Langlands--Shahidi method are in fact Fourier transforms of the measure which defines the corresponding intertwining operator by convolution.
This general point of view was used in [Sh2] to show equality of the Rankin--Selberg $L$--functions and root numbers for $GL_m\times GL_n$ from Langlands--Shahidi method with those defined by Jacquet, Piatetski--Shapiro and Shalika [JPSS].
Whether this Fourier transform can be related to the conjectural ones in [BK,La] remains to be seen.
See Remark 9.4.16.

These notes were presented in part as a series of five lectures at the Morningside Center of Mathematics at Beijing during a workshop organized by Ye Tian and Yangbo Ye in July of 2015.
I would like to thank them for their hospitality and a well--organized workshop.
The last part of these notes, Section 9, was presented during a talk at an AIM workshop organized by Jayce Getz, Dihua Jiang and Lei Zhang over the period November 30--December 4, 2015, for which I like to thank them for the invitation.

\bigskip\noindent
{\bf 1.\ Artin Factors}.
Let $k$ be a number field and let $G_k=\Gal (k_{\alg}/k)$, where $k_{\alg}$ denotes an algebraic closure of $k$.
Let
$$
\rho\colon G_k\longrightarrow GL_n (\bC)=\Aut V,
$$
$V=\bC^n$, be a continuous $n$--dimensional representation.
We can then choose a finite Galois extension $K/k$ such that $\rho$ factors through $G=\Gal(K/k)$.

For each place $v$ of $k$ and $w$ of $K$ with $w|v$, let $\Gal(K_w/k_v)$ be the corresponding decomposition group.
Then the decomposition group will depend, up to conjugation by elements in $G$, only on $v$.
We thus denote it by $G_v$.
Denote by $I_v\subset G_v$ the inertia group at $v$ which again up to conjugation is well--defined and depends only on $v$.
Let $\oK_w$ and $\ok_v$ be the corresponding residue fields with $q_v=\Card (\ok_v)$.
Then
$$
\Gal(\oK_w/\ok_v)\simeq I_v\backslash G_v.
$$
Let $F_v$ be the well defined Frobenius coset of $I_v$, which is again unique up to conjugation.
It is the preimage of a generator of $\Gal(\oK_w/\ok_v)$.

The restriction of $\rho(F_v)$ to $V^{I_v}$, the subspace of $I_v$--invariants in the space $V$ of $\rho$, is independent of $K_w$ and thus
$$
\det (I_V-\rho (F_v)|V^{I_v}\cdot q_v^{-s})^{-1}\tag1.1
$$
is independent of the choice of $K$ and $w$.
Here $s\in\bC$ is a complex variable.
Let $\rho_v$ be the restriction of $\rho$ to $G_v$.
We then define the local Artin $L$--function at $v$ by 
$$
\aligned
L_v(s,\rho)&=L(s,\rho_v)\\
&=\det (I_V-\rho (F_v)|V^{I_v}\cdot q_v^{-s})^{-1}.
\endaligned\tag1.2
$$

For archimedean places, Artin defined an archimedean local $L$--function $L_\infty(s,\rho)$ as a product of $\Gamma$--functions.
We refer to [L4,Sh3,T2] for their definition in the language of automorphic forms.

Next, let $d_{k/\bQ}$ be the discriminant of $k$ over $\bQ$ and let $f=f(\rho)$ be the Artin conductor of $\rho$, an ideal of the ring of integers $O_k$ of $k$ (cf.~[A1,A2]).
It is then a product of local conductors $\prod\limits_{v<\infty} f(\rho_v)$, where $f(\rho_v)$ is the local conductor of $\rho_v$ which depends on the further restrictions of $\rho_v$ to higher ramification groups [A1,A2].
Then Artin defined a global root number $W(\rho)$ which is a root of unity and $\var(s,\rho)$ by
$$
\var(s,\rho)=W(\rho)[|d_{k/\bQ}|^n |N_{k/\bQ} (f(\rho))|]^{-(s-{1\over 2})}\tag1.3
$$
such that the functional equation
$$
L(s,\rho)=\var(s,\rho) L(1-s,\rho^\check),\tag1.4
$$
is valid, where $\rho^\check$ is the dual of $\rho$ (cf.~[A1,A2]).

As it is evident both $d_{k/\bQ}$ and $N_{K/\bQ} (f(\rho))$ are products of local factors.
It was a decomposition of $W(\rho)$ to local factors which resisted efforts until Dwork [Dw], Langlands [L2,L3] and finally Deligne [D] succeeded to define at each place $v$, an $\var$--factor $\var(s,\rho_v,\psi_v)$ such that
$$
\var(s,\rho)=\prod_v\var(s,\rho_v,\psi_v).\tag1.5
$$
We remark that to make the decomposition one needs to fix an additive character $1\neq\psi=\otimes_v \psi_v$ of $k\backslash\bA_k$, where $\bA_k$ is the ring of adeles of $k$.
The choice of $\psi$ is irrelevant for the definition of the global root number, but is needed for its local factors and the decomposition.
This in particular decomposes $W(\rho)$ as
$$
W(\rho)=\prod_v W (\rho_v,\psi_v).\tag1.6
$$
The point is that the local factors are now defined by only local ingredients and have nothing to do with the global data that $\rho_v$ may be obtained from by restriction.

\bigskip\noindent
{\bf 2.\ Local Langlands Correspondence for $GL(n)$ (LLC)}.
Let $F$ be a local field, $p$--adic, real or complex.
Let $W'_F$ be the corresponding Weil--Deligne group [T2].
The local Langlands correspondence (LLC) attaches to every $n$--dimensional continuous complex representation $\rho$ of $W'_F$, for which $\rho$(Frob) acts by semisimple matrices when $F$ is $p$--adic and Frob denotes any element in the inertia coset of a Frobenius element, an irreducible admissible representation $\pi(\rho)$ of $GL_n(F)$ such that the pair $(\rho,\pi(\rho))$ satisfy the following reciprocity rules:

\medskip\noindent
{\it (2.1)\ Identifying abelianization $W_F^{ab}$ of $W_F$ and the center of $GL_n(F)$ with $F^*$, we have
$$
\det (\rho (a))=\omega_{\pi(\rho)}(a)\tag2.2.1
$$
for $a\in F^*$, where $\omega_\pi$ is the central character of $\pi(\rho)$.}

Fix a non--trivial (additive) character $\psi$ of $F$.
Given two irreducible admissible representations $\pi_1$ and $\pi_2$ of $GL_{n_1}(F)$ and $GL_{n_2}(F)$, $n_1,n_2\in\bN$, respectively, let $L(s,\pi_1\times\pi_2)$ and $\var(s,\pi_1\times\pi_2,\psi)$ be Rankin--Selberg (Rankin product) $L$--function and $\var$--factor attached to $\pi_1$ and $\pi_2$ (cf.~[JPSS,Sh2]).
Then

\medskip\noindent
{\it (2.2)\ LLC preserves these $L$--functions and root numbers.
More precisely, if $\rho_1$ and $\rho_2$ are $n_1$--dimensional and $n_2$--dimensional representations of $W'_F$ as discussed before, then
$$
L(s,\rho_1\otimes\rho_2)=L(s,\pi(\rho_1)\times\pi(\rho_2))\tag2.2.1
$$
and
$$
\var (s,\rho_1\otimes\rho_2,\psi)=\var (s,\pi(\rho_1)\times \pi(\rho_2),\psi),\tag2.2.2
$$
where the factors on the left are those of Artin discussed in the previous section.}

LLC has been proved for arbitrary $n$ by Harris and Taylor in [HT].
Soon after that Henniart gave a simpler proof in [He1].
There is now another proof, fairly recent and with a different formulation, due to Peter Scholze [Sch].
All the proofs eventually require a local--global argument to prove (2.2) using functional equations satisfied by global objects which have the local factors in (2.2) attached to their local components.

\medskip\noindent
(2.3)\ {\bf Remark}.
The way LLC is formulated requires the matching of Rankin product $L$--functions and root numbers.
In principle, one may use matching of other $L$--functions and root numbers to formulate the LLC.
But the matching of the Rankin product factors seems to be the most natural and historically the most explored objects for this purpose.
On the other hand with LLC in hand there remains the question of matching other factors under the established LLC.
This is the subject matter of this paper where we explore how the equality of different factors can be established, assuming certain basic axioms for the analytic factors, which are all satisfied by Artin factors, once LLC is established.
We will basically do this for $GL(n)$.
But similar ideas should prevail for arbitrary reductive group.

\bigskip\noindent
{\bf 3.\ An example:\ \ The case of exterior and symmetric square factors for $GL(n)$}

Let $\rho$ be a continuous, $n$--dimensional, Frobenius semisimple complex representation of $W'_F$.
Let $r$ denote either the exterior square $\Lambda^2$ on symmetric square $\Sym^2$ representation of $GL_n(\bC)$.
If $r=\Lambda^2$, then $r|SL_n(\bC)$ will have the second fundamental weight $\delta_2$ of $SL_n(\bC)$ as its highest weight, while for $r=\Sym^2$, the highest weight is $2\delta_1$, where $\delta_1$ is the first fundamental weight.
(More about highest weights later.).
We note that $\dim\Lambda^2={n(n-1)\over 2}$ and $\dim\Sym^2={n(n+1)\over 2}$.
Then $r\cdot\rho$ will be a ${n(n\mp 1)\over 2}$--dimensional representation $W'_F$.
Let $L(s,r\cdot\rho)$ and $\var(s,r\cdot\rho,\psi)$ be the corresponding Artin factors.

In this case there are analytic factors defined by Langlands--Shahidi method [Sh5,Sh8] denoted by $L(s,\pi,r)$ and $\var(s,\pi,r,\psi)$ for every irreducible admissible representation $\pi$ of $GL_n(F)$ which complement the unramified factors of Langlands to provide the functional equations for (global) cuspidal representations which have $\pi$ as a local component.

It is convenient to define
$$
\gamma(s,\pi,r,\psi)=\var(s,\pi,r,\psi)L(1-s,\tilde\pi,r)/L(s,\pi,r).\tag3.1
$$
Similarly we set
$$
\align
& \gamma(s,r\cdot\rho,\psi)=\var (s,r\cdot\rho,\psi)L(1-s,r^\check\cdot\rho)/L(s,r\cdot\rho),\tag3.2\\
& \gamma(s,\pi_1\times\pi_2,\psi)=\var(s,\pi_1\times\pi_2,\psi)
L(1-s,\tilde\pi_1\times\tilde\pi_2)/L(s,\pi_1\times\pi_2),\tag3.3
\endalign
$$
as well as
$$
\gamma(s,\rho_1\otimes\rho_2,\psi)=\var(s,\rho_1\otimes\rho_2,\psi)L(1-s,\rho_1^\check\otimes\rho_2^\check)/L(s,\rho_1\otimes\rho_2).\tag3.4
$$
Here $\rho^\check$ is the dual of $\rho$.

The following is proved in [CST]:

\NoBlackBoxes
\medskip\noindent
(3.5)\ {\bf Theorem} (Cogdell--Shahidi--Tsai [CST]).
{\it Let $\rho$ be an $n$--dimensional continuous Frobenius--semisimple complex representation of $W'_F$ and assume $r=\Lambda^2$ or $r=\Sym^2$.
Then
$$
\var(s,r\cdot\rho,\psi)=\var(s,\pi(\rho),r,\psi)\tag3.5.1
$$
and
$$
L(s,r\cdot\rho)=L(s,\pi(\rho),r).\tag3.5.2
$$}

\medskip\noindent
(3.6)\ {\bf Remarks}:\ \ a) The equality (3.5.2) was first proved by Henniart [He2].
The case of root numbers (3.5.1) is much harder and requires stability of root numbers to be discussed soon.

b)\ Theorem 3.5 is proved by first establishing 
$$
\gamma (s,r\cdot\rho,\psi)=\gamma(s,\pi(\rho),r,\psi)\tag3.6.1
$$
and deducing (3.5.1) and (3.5.2) by appealing to the machinary developed in [Sh5] through Langlands classification [L4,Si1,Si2].

c)\ One can show
$$
\gamma(s,\pi\times\pi,\psi)=\gamma(s,\pi,\Lambda^2,\psi)\gamma(s,\pi,\Sym^2,\psi)\tag3.6.2
$$
as well as
$$
\gamma(s,\rho\times\rho,\psi)=\gamma(s,\Lambda^2\cdot\rho,\psi)\gamma(s,\Sym^2\cdot\rho,\psi).\tag3.6.3
$$
Consequently one only needs to prove (3.6.1) for $r=\Lambda^2$ by (2.2.1) and (2.2.2) which implies
$$
\gamma(s,\rho\otimes\rho,\psi)=\gamma(s,\pi(\rho)\times\pi(\rho),\psi).\tag3.6.4
$$

\bigskip\noindent
{\bf 4.\ The Proof}.
In this section we will outline how the proof of Theorem 3.5 proceeds.
The technology that goes into the proof is very robust and can be generalized to other representations of $GL_n(\bC)$, whenever an analytic theory whose $\gamma$--factors satisfy a number of properties as those used here, exists.
This will be the subject matter of the next section.
We now start with different properties and tools needed and used for the proof.

\medskip\noindent
(4.1)\ {\bf Additivity/Multiplicativity}.
Arithmetic (Artin) factors satisfy a number of properties in full generality.
In the case where $r=\Lambda^2$ one has
$$
\gamma(s,\Lambda^2 (\rho_1\oplus\rho_2),\psi)=\gamma(s,\Lambda^2\rho_1,\psi)\gamma(s,\Lambda^2\rho_2,\psi)\gamma(s,\rho_1\otimes\rho_2,\psi),\tag4.1.1
$$
where $\rho_i,\ i=1,2$, are continuous representations of $W'_F$ as before.
This is the additivity for these factors.
We refer to [L2,L3,CST] for this and other properties of Artin factors.

On the analytic side this becomes
$$
\gamma(s,\Ind(\pi_1\otimes\pi_2),\Lambda^2,\psi)=\gamma(s,\pi_1,\Lambda^2,\psi)\gamma(s,\pi_2,\Lambda^2,\psi)\gamma(s,\pi_1\times\pi_2,\psi).\tag4.1.2
$$
This is what is usually called ``multiplicativity'' and is proved for all the factors defined by Langlands--Shahidi method in [Sh5].
We will denote this property by $(M)$.

\medskip\noindent
(4.2)\ {\bf Stability}.
One of the tools used by Deligne [D] to define the local $\var$--factors and root numbers was stability of $\var$--factors under twisting by a highly ramified characters, i.e., a character $\chi$ of $F^*$ whose conductor is suitably large.

\medskip\noindent
{\bf Arithmetic Stability} (Deligne [D]).
Let $\rho$ be a representation of $W'_F$ as before.
Then there exists a positive integer $f_0$, depending on $\rho$, such that for each character $\chi$ of $F^*\simeq W_F^{ab}$ whose conductor $f\geq f_0$, there exists a $y\in F^*$, depending on $\chi$ and $\psi$, such that 
$$
\var(s,\rho,\psi)=(\det\rho)^{-1} (y)\var (s,\chi,\psi)^{\dim\rho}.\tag4.2.1
$$
In particular, if $\rho_1$ and $\rho_2$ are two representations of $W'_F$ as before with $\det\rho_1=\det\rho_2$, then for all highly ramified characters $\chi$ of $F^*$
$$
\var(s,\rho_1\otimes\chi,\psi)=\var(s,\rho_2\otimes\chi,\psi)\tag4.2.2
$$
with the degree of ramification $f_0$ depending on $\rho_1$ and $\rho_2$.

In the case of exterior square for $GL_n$ we can write (4.2.2) in terms of $\gamma$--functions and in the form:
$$
\gamma(s,\Lambda^2 (\rho_1\otimes\chi),\psi)=\gamma (s,\Lambda^2 (\rho_2\otimes\chi),\psi),\tag4.2.3
$$
whenever $\det\rho_1=\det\rho_2$ and $\chi$ is highly ramified.

\medskip\noindent
{\bf Analytic Stability} [CPS,CPSS1,CPSS2,JS,Sh6,CST].
A similar situation is expected to happen in the analytic side as well.
Let $\pi_1$ and $\pi_2$ be two irreducible admissible representations of $GL_n(F)$ with equal central characters $\omega_{\pi_1}=\omega_{\pi_2}$.
Let $r$ be a finite dimensional representation of $GL_n(\bC)$ and assume there are defined analytic factors $\gamma(s,\pi_i,r,\psi),i=1,2$.
Then stability demands that for all highly ramified characters $\chi$ of $F^*$
$$
\gamma(s,\pi_1\otimes\chi,r,\psi)=\gamma(s,\pi_2\otimes\chi,r,\psi).\tag4.2.4
$$
Again the degree of ramification will depend on $\pi_1$ and $\pi_2$.

When $\pi_i$ are supercuspidal (4.2.4) is called ``Supercuspidal Stability'' and is denoted by {\bf (SCS)}.
As will be explained later it is much easier to prove SCS, at least in examples.
Moreover, the method which will be proposed here to prove (3.6.1) for any given $r$, as long as these factors are defined satisfying a number of properties, will only require (SCS) as far as stability is concerned.
Stability for arbitrary irreducible admissible representations $\pi_1$ and $\pi_2$ then follows from (3.6.1) for the given $r$, using arithmetic stability.

When $r=\Lambda^2$, stability or (SCS) will be 
$$
\gamma(s,\pi_1\otimes\chi,\Lambda^2,\psi)=\gamma(s,\pi_2\otimes\chi,\Lambda^2,\psi),\tag4.2.5
$$
whenever $\omega_{\pi_1}=\omega_{\pi_2}$ and $\chi$ is highly ramified.

\medskip\noindent
(4.3)\ {\bf Archimedean Matching (AM)}.
The next ingredient needed in the proof of Theorem 3.5 is the equality (3.6.1) when $F$ is an archimedean field and $r=\Lambda^2$.
But this is a theorem proved in full generality when $r$ is one of the representations of the $L$--group appearing in the Langlands--Shahidi method, the case of $r=\Lambda^2$ being among them.
The equality (3.6.1) when $F$ is archimedean is called archimedean matching and is denoted by {\bf (AM)}.

\medskip\noindent
(4.4)\ {\bf Functional Equation} {\bf (FE)}.
This is the only global tool which goes to the proof of Theorem 3.5 and is repeatedly used in different steps.
Here one needs it both for global representations of Weil group and cuspidal automorphic forms for $GL_n(\bA_k)$, where $k$ is a number field with $\bA_k$ its ring of adeles.
If $\rho$ is a continuous $n$--dimensional representation of $G_k$ as in Section 1, then the functional equation satisfied by $\rho$ is equation (1.4) of Section 1 with factors as defined there.

Now let $\pi=\otimes_v\pi_v$ be a cuspidal automorphic representation of $GL_n(\bA_k)$, where each $\pi_v$ is an irreducible unitary representation of $GL_n(k_v)$.
Since $\pi$ has a non--vanishing Fourier coefficient with respect to any generic character of $U(k)\backslash U(\bA_k)$, where $U$ is the subgroup of upper triangular unipotent matrices in $GL_n$, then each $\pi_v$ is generic ([Sha]).
Now let $r$ be a finite dimensional representation of $GL_n(\bC)$.
Let $S$ be a finite set of places of $k$ such that for each $v\not\in S$, $\pi_v$ is unramified, i.e., the space $\sH(\pi_v)$ contains a vector fixed by $GL_n(O_v)$ with $O_v$ the ring of integers of $k_v$.
Such representations are then parametrized by semisimple conjugacy classes in $GL_n(\bC)$.
Let $t_v$ be the one attached to $\pi_v$.
Langlands defined an $L$--function attached to $\pi_v$ and $r$ by
$$
L(s,\pi_v,r)=\det(I-r(t_v) q_v^{-s})^{-1}\tag4.4.1
$$
which clearly depends only on the conjugacy class $t_v$ and not any representative of it.
Now, define
$$
L^S(s,\pi,r)=\prod_{v\not\in S} L(s,\pi_v,r).\tag 4.4.2
$$
One then expects for each $v$ the existence of scalar functions $\gamma(s,\pi_v,r,\psi_v)$, with $\psi_v$ a local component of a global non--trivial character $\psi=\otimes_v\psi_v$ of $k\backslash \bA_k$, such that
$$
L^S(s,\pi,r)=\prod_{v\in S}\gamma (s,\pi_v,r,\psi_v) L^S(1-s,\tilde\pi,r),\tag4.4.3
$$
where $\tilde\pi$ is the contragredient of $\pi$.
We point out that when $\pi_v$ is unramified
$$
\gamma(s,\pi_v,r,\psi_v)=L(1-s,\tilde\pi_v,r)/L(s,\pi_v,r),\tag4.4.4
$$
which in our context implies $\var(s,\pi_v,r,\psi_v)=1$ for any unramified $\pi_v$.
This is the functional equation satisfied by $L^S(s,\pi,r)$ and will be denoted by (FE).
We note that again when $r$ is one of the representations of $GL_n(\bC)$ appearing in the Langlands--Shahidi method (or more generally of $^LM$, cf.~[L1,Sh8]), then (4.4.3) is always valid and a theorem proved [Sh1,Sh5], referred to as (FE).

\bigskip\noindent
5.\ {\bf Steps of the proof}.
Our proof will reduce the proof of Theorem 3.5, which is basically a proof of (3.6.1) for $r=\Lambda^2$, to a proof of (SCS).
We will now sketch these steps and later that of (SCS) in this case.

We will first establish a stable version of (3.6.1), namely:

\medskip\noindent
(5.1)\ {\bf Proposition}.
{\it Fix $n\in\bN$ and let $\rho$ be an $n$--dimensional irreducible complex representation of $W_F$.
Then for all highly ramified characters $\chi\in\hat F^*$, one has}
$$
\gamma(s,\Lambda^2(\rho\otimes\chi),\psi)=\gamma(s,\pi(\rho)\otimes\chi,\Lambda^2,\psi).\tag5.1.1
$$

To prove the proposition, we first show the existence a base point $(\rho_0,\pi(\rho_0))$ for which (5.1.1) is valid for all $\chi\in\hat F^*$.
This will then allow us a deformation to prove (5.1.1) for any $(\rho,\pi(\rho))$.
We need

\medskip\noindent
(5.2)\ {\bf Lemma}.
{\it Given $n\in\bN$ and $\omega_0\in\hat F^*$, there exists a number field $k$, an $n$--dimensional complex representation $\tilde\rho$ of $W_k$ such that, if for each place $w$ of $k$, $\tilde\rho_w$ denotes the restriction of $\tilde\rho$ to $W_{k_w}$, then}

\medskip\noindent
(5.2.1)\ {\it there exists a finite place $v$ of $k$ such that $k_v=F$, $\tilde\rho_v$ is irreducible and $\det\tilde\rho_v=\omega_0$;}

\medskip\noindent
(5.2.2)\ {\it for every finite place $w\neq v$, $\tilde\rho_w$ is reducible, and}

\medskip\noindent
(5.2.3)\ {\it the irreducible admissible representation $\pi(\tilde\rho)\colon = \otimes_w\pi(\tilde\rho_w)$ is a cuspidal automorphic representation of $GL_n(\bA_k)$.}

The proof of this lemma was communicated to us by Henniart upon our inquiry to its validity.

With Lemma 5.2 in hand we will now compare the functional equations for $\tilde\rho\otimes\tilde\chi$ and $\pi(\tilde\rho\otimes\tilde\chi)$, where $\tilde\chi=\otimes_w \tilde\chi_w$ is a gr\"ossencharacter of $k$, which is highly ramified at all finite ramified places different from $v$, while $\tilde\chi_v=\chi$ with $\chi$ as in Proposition 5.1.

We now use induction by assuming the validity of Proposition 5.1 for any local field $F$ and any $m<n$.
We can then use equation (4.1.2), multiplicativity, to show that
$$
\gamma(s,\Lambda^2(\tilde\rho_w\otimes\tilde\chi_w),\psi_w)=\gamma(s,\pi(\tilde\rho_w\otimes\tilde\chi_w),\Lambda^2,\psi_w)\tag5.2.4
$$
for all $w\neq v$.
The validity of (5.2.4) for all $w=\infty$ and any $\tilde\chi_w$ follows from archimedean matching (AM) in (4.3).
We now compare functional equations for $\tilde\rho\otimes\tilde\chi$ and $\pi(\tilde\rho\otimes\tilde\chi)$ to conclude:

\medskip\noindent
(5.3)\ {\bf Conclusion} {\it (Existence of a base point for deformation).
There exists a pair \newline
$(\rho_0,\pi(\rho_0))$ with $\rho_0$ irreducible and thus $\pi_0$ supercuspidal, such that
$$
\gamma(s,\Lambda^2 (\rho_0\otimes\chi),\psi)=\gamma(s,\pi(\rho_0)\otimes\chi,\Lambda^2,\psi)\tag5.3.1
$$
for all $\chi\in\hat F^*$}.

We now assume $\chi$ is highly ramified.
Then by arithmetic stability (4.2.3) we have that
$$
\gamma(s,\Lambda^2 (\rho\otimes\chi),\psi)=\gamma (s,\Lambda^2 (\rho_0\otimes\chi),\psi)\tag5.3.2
$$
for all highly ramified $\chi$, with ramification depending on $\rho$ and $\rho_0$.

Next assume (SCS) for $\rho$ to the effect of the validity of (4.2.5) which we display as
$$
\gamma(s,\pi(\rho)\otimes\chi,\Lambda^2,\psi)=\gamma(s,\pi(\rho_0)\otimes\chi,\Lambda^2,\psi)\tag5.3.3
$$
Putting (5.3.1), (5.3.2) and (5.3.3) together we now have
$$
\gamma(s,\Lambda^2(\rho\otimes\chi),\psi)=\gamma(s,\pi(\rho)\otimes\chi,\Lambda^2,\psi)\tag5.3.4
$$
for any irreducible $\rho$ and with $\chi$ highly ramified depending on it.
This proves Proposition 5.1, but subject to validity of (SCS) for $r=\Lambda^2$, which we will sketch a proof for soon.

It follows from arithmetic and analytic multiplicativity that

\medskip\noindent
(5.4)\ {\bf Corollary}.
{\it Proposition 5.1 is valid for any $n$--dimensional continuous Frobenius semisimple representation $\rho$ of $W'_F$.}

With Corollary 5.4 in hand we can now sketch the proof of Theorem 3.5.
We will first prove the theorem for a $\bZ$--basis of the Grothendieck group of all finite dimensional representations of $W_F$ by Brauer's theorem [Br], the monomial representations, i.e., those of the form $\rho=\Ind_H^G\eta$, where $\eta$ is a character of a subgroup $H$ of $G$ of finite index.
Following Harris (cf.~[H,He2,CST]), there exist a number field $k$ and a finite dimensional representation $\tilde\rho$ of $W_k$ such that there exists a place $v$ of $k$ for which $k_v=F$ and $\tilde\rho|k_v=\rho$, is a fixed monomial representation of $W_F$.
Moreover, $\tilde\pi\colon=\otimes_w \pi (\tilde\rho_w)$ is an automorphic representation of $GL_n(\bA_k)$.
We then let $\tilde\chi=\otimes_w \tilde\chi_w$, where $\tilde\chi_v=1$, while $\tilde\chi_w$ is highly ramified for every $w\neq v,\ w<\infty$, for which $\tilde\pi_w=\pi(\tilde\rho_w)$ is ramified.
We then compare functional equations for $L(s,\Lambda^2(\tilde\rho\otimes\tilde\chi))$ with $L(s,\tilde\pi\otimes\tilde\chi,\Lambda^2)$, using Corollary 5.4.
We note that when $k_w$ is an archimedean field, then (5.3.4) is valid for any $\rho$ and $\chi$ by the general results proved in [Sh3], using the LLC in [L4].
Equality for any monomial representation then follows.

Proof for an arbitrary $\rho$ now follows from the equality
$$
\Lambda^2 (\rho_1\oplus\rho_2)=\Lambda^2 \rho_1\oplus\Lambda^2 \rho_2\oplus\rho_1\otimes\rho_2\tag5.4.1
$$
coming from additivity/multiplicativity from which one concludes
$$
\Lambda^2 (\rho_1\ominus\rho_2)=\Lambda^2\rho_1\ominus\Lambda^2\rho_2\ominus\rho_1\otimes\rho_2\oplus\rho_2\otimes\rho_2.\tag5.4.2
$$
Our theorem now follows if $\rho_1$ and $\rho_2$ are monomial representations of any given dimension.

\medskip\noindent
(5.5)\ {\bf Remark}.
We note that to prove the theorem for arbitrary $\rho$ of dimension $n$ one needs its validity for all the monomial representation, and not only those of dimension $n$.

\medskip\noindent
(5.6)\ {\bf Corollary}.
{\it Analytic stability is valid for all irreducible admissible representations of $GL_n(F)$ when $r=\Lambda^2$ or $\Sym^2$}.

This follows from stability of arithmetic factors as proved by Deligne for any $r$ (cf.~Section 4.2).

\bigskip\noindent
6.\ {\bf Proof of} {\bf (SCS)}.
We now sketch the proof of analytic stability for supercuspidal representations (SCS) given in [CST].
The proof of this is of different nature and rather hard.
We need to study the asymptotic behavior of certain partial Bessel functions [Sh6,CST].
While (full) Bessel functions have germ expansion due to Jacquet--Ye [J1,JY], the partial ones lack one, and that is where difficulty lies.
Fortunately, we could use a number of ideas in [J1,JY] to prove an asymptotic expansion, alas not a germ expansion, for our partial Bessel function to prove stability in this case.
As far as we know, no germ expansion is expected for our partial Bessel functions.

Analytic factors coming from Langlands--Shahidi method in general and in particular in this case are defined by local coefficients attached to a triple $(H,M_H,\sigma)$, where $H$ is a quasisplit connected reductive group over $F$ with $M_H$ a Levi subgroup of a maximal parabolic subgroup of $H$, again over $F$, and $\sigma$ is an irreducible generic representation of $M_H(F)$.
We refer to [L1,Sh1,Sh5,Sh8] for the details on the defining objects.

In the case of $\gamma(s,\pi,\Lambda^2,\psi)$, we can take $H=GSp_{2n}$, $M_H=GL_n\times GL_1$, the Siegel Levi subgroup, $\sigma=\pi_s\otimes\bold 1$, where $\pi_s=\pi\otimes|\det(\quad)|^s$, $s\in\bC$.
There are other choices such as $H=SO_{2n}$ or $H=Sp_{2n}$ and $M_H=GL_n$.
But when $n$ is odd this Levi is not self--associate in $SO_{2n}$ and the center of $Sp_{2n}$ is not connected, both conditions which must be satisfied if we are to use the main integral formula proved in [Sh6], Theorem 6.2, to prove stability.

The defining local coefficient $C_\psi(s,\pi)$ for this triple equals
$$
C_\psi (s,\pi)=\gamma (2s,\tilde\pi,\Lambda^2,\overline\psi)\gamma(s,\tilde\pi,\stand,\opsi),\tag6.1
$$
where the $\gamma$--function $\gamma(s,\pi$, stand, $\psi)$ is the Godement--Jacquet $\gamma$--factor attached to the standard representation of $GL_n(\bC)$, the $L$--group of $GL_n$, in [GJ].

Theorem 6.2 of [Sh6] expresses $C_\psi (s,\pi\otimes\chi)^{-1}$, up to an abelian $\gamma$--factor depending only on $\omega_\pi$ and $\chi$, as a Mellin transform of a partial Bessel function, whenever $\chi$ is sufficiently ramified as we explain later.
For simplicity we use PBF and MT do denote the partial Bessel function and its Mellin transform, letting us to write 
$$
C_\psi (s,\pi\otimes\chi)^{-1}\sim MT (PBF)\tag6.2
$$
for brevity.
The supercuspidality of $\pi$ now plays a central role since its matrix coefficients are of compact support modulo the center $F^*$ of $GL_n(F)$, allowing us to pursue the ideas in [J1,JY,Sh6] to prove (SCS).
On the other hand the case of general $\pi$ seems to be much harder since the asymptotics of PBF in general seems very much out of reach.

\medskip\noindent
(6.3)\ {\bf Partial Bessel functions}.
Let $G=GL_n(F)$ and let $A$ be the subgroup of diagonals and $B=AU$ the Borel subgroup of upper triangulars with $U$ its unipotent radical.
Let $Z$ be the center of $G$ and
$$
A'=\left\{\pmatrix 1\\
& a_2 & & 0\\
&&\ddots\\
0&&&a_n
\endpmatrix \Bigg| a_i\in F^*\right\}\subset A.\tag6.3.1
$$

Assume $\pi$ is irreducible supercuspidal and let $\omega$ be a character of $F^*$ which most of the time will be $\omega_\pi$, the central character of $\pi$.
Denote by $C_c^\infty(G,\omega)$ the space of smooth functions of compact support modulo center, transforming according to $\omega$ through identification $Z\simeq F^*$.
When $\omega=\omega_\pi$ we can identify the space of matrix coefficient $\sM(\pi)$ of $\pi$ as a subspace $\sM(\pi)\subset C_c^\infty(G,\omega)$, since $\pi$ is supercuspidal.

Let $w\in W(G,A)$ which by abuse of notation denotes the Weyl group of $A$.
Given $w$, let $U_w^-$ be the subgroup of $U$ generated by those simple root vectors which under the action of $w$ go to a negative roots.
Let $C(w)=U\dot w A U_w^-$, be the corresponding Bruhat cell (double coset), where $\dot w$ is a representative for $w$ fixed as in [Sh5,Sh8].
Given $g\in C(w)$, write $g=zg'$, $z\in Z$, $g'\in C'(w)=U\dot w A' U^-_w$.
The decomposition is then unique.

Given $u\in U$, we use $\psi$ to denote the generic character of $U$ defined by
$$
\psi(u)=\psi(u_{1,2} + u_{2,3} + \ldots + u_{n-1,n}).\tag6.3.2
$$
Fix $f\in\sM (\pi)$ and for $g\in C(w)$, define
$$
W^f(g)=\int_U f(xg)\psi^{-1} (x) dx.\tag6.3.3
$$

Recall that $H=GSp_{2n}$ and $M_H=GL_n\times GL_1$.
Let $P_H=M_H N_H$ be the standard Siegel parabolic subgroup of $H$ with unipotent radical $N_H$.
We let $\oN_H$ be its opposite unipotent subgroup.
Then
$$
\oN_H=\left\{ \overline n(y)=\pmatrix I_n&0\\ y&I_n\endpmatrix\Bigg| y=J^t yJ\right\},\tag6.3.4
$$
where $J=\pmatrix
&&&1\\
0&&-1\\
&\ddots\\
(-)^{n-1}&&0
 \endpmatrix$.
We recall from [CST] the choices made in defining $H$ as
$$
H(F)=GSp_{2n} (F)=\{h\in GL_{2n} (F)|^th J' h=\eta(h) J',\exists \eta (h)\in F^*\},\tag6.3.5
$$
$$
J=J_n=\pmatrix
&&&1\\
&&-1\\
&\ddots\\
(-1)^{n-1}
\endpmatrix\tag6.3.6
$$
and
$$
J'=J'_{2n}=\pmatrix 0&J\\ -{}^t\!J&0\endpmatrix.\tag6.3.7
$$
Finally we have our representative $\dot w_\ell$ for the long element of $W(G,A)$ as in [Sh8].
It equals $\dot w_\ell=J$.

Fix a positive integer $N$.
We now define a cutoff function $\varphi_N$ on $M_n(F)$ by
$$
\varphi_N=\char X(N),\tag6.3.8
$$
where
$$
X(N)=\left\{\pmatrix\ff^{-N}&\ff^{-2N}&\ff^{-3N}\\
\ff^{-2N}&\ff^{-3N}&\ff^{-4N}&\ldots\\
\ff^{-3N}&\ff^{-4N}&\ff^{-5N}\\
&\ddots
\endpmatrix\right\}.\tag6.3.9
$$
with the prime ideal $\ff$ of the ring of integers $O$ of $F$.
We note that the sets $\tilde X(N)\cap N_H(F)$ exhaust $N_H(F)$ as $N$ increases, where
$$
\tilde X(N)=\left\{ \pmatrix I_n& X(N)\\ 0&I_n\endpmatrix\right\} \subset M_{2n} (F).\tag6.3.10
$$

Given $g\in G$, let
$$
\aligned
U_g&=\{ u\in U| ^tu \dot w_{\ell}^{-1} gu=\dot w_{\ell}^{-1} g\}\\
&= \left\{ u\in U|\tilde u g u=g,\tilde u=\dot w_{\ell} ^tu \dot w_{\ell}^{-1} \right\}
\endaligned\tag6.3.11
$$
be the twisted centralizer of $g$.
Then $U_{\dot w_\ell}=U_{\dot w_\ell a}= \{e\}$, for all $a\in A$.
Let
$$
U(N)=\{ u=(u_{ij})\in U| u_{ij}\in \ff^{-N},\ i<j,\ i\geq 1 \}.\tag6.3.12
$$
Then we have
$$
\varphi_N (^tu gu)=\varphi_N (g),\tag6.3.13
$$
for all $u\in U(N)$.
We can now define the {\it partial Bessel function} attached to $f\in\sM(\pi)$ and $\varphi=\varphi_N$ for some $N\in\bN$ as 
$$
\aligned
B_\varphi^G(g,f)&=\int_{U_g\backslash U} W^f (gu) \varphi (^t u\dot w_{\ell} g' u)\psi^{-1} (u) du\\ \\
&=\int_{U_g\backslash U} \int_U f(xgu)\varphi(^tu\dot w_{\ell}^{-1} g' u)\psi^{-1} (xu) dx du\\ \\
&=\omega_\pi (z) B_\varphi^G (g',f),
\endaligned\tag6.3.14
$$
where $g\in U\dot w A U_w^-,\ g=zg',\ g'\in U\dot w A' U_w^-$.
Then for $a\in A$, we have
$$
B^G_\varphi (\dot w_\ell a,f)=\int_U W^f (\dot w_\ell a u)\varphi (^tu a' u)\psi^{-1} (u) du.\tag6.3.15
$$
We note that
$$
B^G_{\varphi_N} (u_1 g u_2,f)=\psi(u_1) B^G_{\varphi_N} (g,f) \psi (u_2)\tag6.3.16
$$
for all $u_1\in U$ and $u_2\in U(N)$, justifying the name ``partial'' Bessel function.

\bigskip\noindent
(6.4)\ {\bf The Mellin transform}.
Theorem 6.2 of [CST] implies that in the present case
$$
\gathered
C_\psi(s,\pi\otimes\chi)^{-1}=\gamma ({ns\over 2},\omega_\pi\chi^n,\psi)\omega_\pi\chi^n (-1)\\
\cdot\int_{Z\backslash A} B_\varphi^G (\dot w_\ell a, f\otimes\chi) \omega_\pi\chi^n (a_1)^{-1} \mu_s (a') da',
\endgathered\tag6.4.1
$$
where $\mu_s$ is a character of $A'$ depending on $s$.
Throughout $f\in\sM(\pi)$ is assumed to satisfy $W^f(e)=1$.

The integral in (6.4.1) then can be written as 
$$
\int_{A'} B_\varphi^G (\dot w_\ell a', f\otimes\chi) \mu_s (a') da'.\tag6.4.2
$$
(The factor $\gamma ({ns\over 2},\omega_\pi\chi^n,\psi)$ is the Hecke--Tate $\gamma$--function attached to $\omega_\pi \chi^n$.)

The integral (6.4.2) is then equal 
$$
\int_{A'} B^G_\varphi (\dot w_\ell a', f)\chi (\det a') \mu_s (a') da',\tag6.4.3
$$
which is clearly the Mellin transform of the partial Bessel function $B^G_\varphi (\dot w_\ell a', f)$.

\medskip\noindent
(6.4.4)\ {\bf Remark}.
The Mellin transform is an integration over the ``big cell'' $B w_\ell B$ modulo the center and to analyze it to prove (SCS) one needs to understand the asymptotic behavior of $B_\varphi^G (-,f)$ as we approach the boundary of the big cell, i.e., every other Bruhat cell.
This is a much more complicated situation than those which appeared in the cases of functoriality and the reason for establishing (SCS) first.
This asymptotics can be expressed, without going to any details as:

\noindent
{\it Up to a non--smooth term which depends on $\pi$ only through $\omega_\pi$ and for a $\varphi$ with support sufficiently large depending on $f$, $B^G_\varphi (\dot w_\ell a,f)$ can be written as a finite sum $\sum\limits_i B^G_\varphi (\dot w_\ell a,f_i)$, where each term is ``uniformly smooth'' on the quotients $Z\backslash A_i$ of certain subtori $A_i\subset A$ by $Z$.}

A function $F$ on $A$ is {\it uniformly smooth} in our context if there exists a subtorus $T, Z\subsetneq T\subset A$ and an open compact subgroup $O$ of identity in $Z\backslash T$ such that $F(at)=F(a)$ for all $a\in A$ and $t\in O$.
This is particularly the case if $F|T$ has a compact support modulo $Z$.

The Mellin transform of $F$ in our context is an integral of the form
$$
MT(F)=\int_{Z\backslash A} F(a) \chi(a) |a|^s d^* a,\tag6.4.6
$$
where $s\in\bC$ and $\chi$ is a character of $A$ and $|a|^s$ is meant to be $|\det a|^s$.
Now fix $t\in O$ and assume $\chi$ is sufficiently ramified such that $\chi(t)\neq 1$ for some $t\in O$.
The subgroup $O$ being compact in $Z\backslash T$ implies $|t|=1$.
A change of variable $a\mapsto at$ in (6.4.6) now implies that $MT(F)=\chi(t) MT(F)$ which implies $MT(F)=0$ since $\chi(t)\neq 1$.

We now explain how this discussion when applied to differences of partial Bessel functions for different representations can lead to a proof of (SCS).

Let $\pi$ and $\pi'$ be two irreducible supercuspidal representations of $G=GL_n(F)$ such that $\omega_\pi=\omega_{\pi'}=\omega$.
Let $f\in\sM(\pi)$ and $f'\in\sM(\pi')$ be matrix coefficients such that 
$$
W^f (e)=W^{f'} (e)=1.\tag6.4.7
$$

Using (6.4.5) and the above discussion on Mellin transforms of uniformly smooth functions applied to (6.4.3), now implies that
$$
MT(B^G_\varphi (-, f\otimes\chi-f'\otimes\chi))=0\tag6.4.8
$$
for a highly ramified character $\chi$ if the support of $\varphi$ is suitably large, depending on $f$ and $f'$, such that (6.4.5) is valid.

It thus follows from (6.1), (6.2) or more precisely (6.4.1) that
$$
\aligned
\gamma(2s,\pi\otimes\chi,\Lambda^2,\psi)\gamma(s,\pi\otimes\chi,\stand,\psi)&=\gamma(2s,\pi'\otimes\chi,\Lambda^2,\psi)\\
&\cdot\gamma(s,\pi'\otimes\chi,\stand,\psi).
\endaligned\tag6.4.9
$$
The equality of the factors attached to the standard representation of $GL_n(\bC)$, the second factors in (6.4.9), is already proved by Jacquet and Shalika in [JS].
This thus gives (SCS) for $\gamma(s,\pi,\Lambda^2,\psi)$ completing the proof of Theorem 3.5.
We discuss the validity of (6.4.5) in the next subsection.

\noindent
(6.5)\ {\bf Asymptotic expansion for partial Bessel functions}.
We now elaborate on (6.4.5) and how the uniform smoothness is proved for our PBF.

We will first address the notion of a Weyl group element $w\in W(G,A)$ supporting Bessel functions.

\medskip\noindent
(6.5.1)\ {\bf Definition}.
{\it An element $w\in W(G,A)$ is said to support a Bessel function if $w_\ell w$ is the long element in the Weyl group of a standard Levi subgroup of $G$ containing $A$.
We denote the set of all such elements by $B(G)$.
(More details to come.)}

This is an important property since full Bessel functions are only supported on Bruhat cells attached to such elements.
The notion plays an important role even for partial Bessel functions as shown in [CST].

For the sake of exposition we will first present the case of $n=3$ and consider $G=GL_3(F)$.
The only elements of $B(G)$ in this case are $w_e,w_c, w_d$ and the long element $w_\ell$, where 
$$
w_c=w_d^{-1}=\pmatrix 0&0&1\\ 1&0&0\\ 0&1&0\endpmatrix.\tag6.5.2
$$
We note that for $w=w_c$ or $w_d$, $w_\ell w$ will be the long elements of the two standard maximal Levi subgroups of $GL_3(F)$, namely $M_c=GL_2 (F)\times GL_1(F)$ and $M_d=GL_1 (F)\times GL_2(F)$, respectively.
We will then denote the centers of $M_c$ and $M_d$ by $A_c$ and $A_d$, respectively.
Then 
$$
A_c=\left\{\pmatrix a\\ &a\\ &&b\endpmatrix\Big| a,b\in F^*\right\}\tag6.5.3
$$
and
$$
A_d=\left\{\pmatrix a\\ &b\\ &&b\endpmatrix\Big| a,b\in F^*\right\}.\tag6.5.4
$$
We note that $A$ is then the center of the Levi subgroup of the standard minimal parabolic subgroup of $G$ and the one attached to $w=w_\ell$.

Given $w\in W(G,A)$, define
$$
\Omega_w=\bigcup_{w\leq w'} C(w'),\tag6.5.5
$$
where $w\leq w'$ is the {\it Bruhat order} on $W(G,A)$, i.e., $w\leq w'$ if and only if $C(w)\subset\overline{C(w')}$.
Every Bruhat cell $C(w)$ lies in $G=\overline{C(w_\ell)}$ and thus $w\leq w_\ell$ for all $w\in W(G,A)$.
Finally, note that $\Omega_w$ is open in $G$ (cf.~[J1,JY]).

A precise version of (6.4.5) in this case is 

\medskip\noindent
(6.5.5)\ {\bf Proposition} {\it (Asymptotic expansions for partial Bessel functions on $GL_3(F)$):\ \ Given $f\in C_c^\infty (GL_3 (F),\omega_\pi)$ with $W^f(e)=1$ and with support of $\varphi$ sufficiently large depending on $f$, there exist functions $f_1\in C_c^\infty (GL_3 (F),\omega_\pi)$, $f_{1,c}\in C_c^\infty (\Omega_{w_c}, \omega_\pi)$, $f_{1,d}\in C_c^\infty (\Omega_{w_d}, \omega_\pi)$ and $f_2\in C_c^\infty (\Omega_{w_\ell},\omega_\pi)$ such that,
$$
\aligned
B^G_\varphi (\dot w_\ell a,f)&=B^G_\varphi (\dot w_\ell a,f_1)+B^G_\varphi (\dot w_\ell a,f_{1,c})\\
&+B^G_\varphi (\dot w_\ell a, f_{1,d})+B^G_\varphi (\dot w_\ell a,f_2).
\endaligned
$$
The first term depends only on $\omega_\pi$ while the last three are uniformly smooth on $Z\backslash A_c$, $Z\backslash A_d$ and $Z\backslash A$, respectively.}

We note that $C(w)$ is closed in $\Omega_w$ for each $w$ which gives the necessary uniformity.

The general case is similar.
The sum will be over all $w$ which support Bessel functions, i.e., whenever $w=w_\ell w^M_\ell$, where $w^M_\ell$ is the long element of a Levi subgroup $M$ of $G$, $M\supset A$, for which simple roots of $M$ for $A$ are among those of $A$ in $U$.
We will call such a Levi subgroup, a {\it standard} one, or a Levi subgroup of a standard parabolic subgroup of $G$ containing $A$.
The elements $w_\ell^M$ are called {\it relevant} and their collection is denoted by $R(G)$ in [J1,JY].
Notice that $R(G)=w_\ell B(G)$.

As for representatives, we must have
$$
\dot w=\dot w_\ell (\dot w_\ell^M)^{-1}\tag6.5.6
$$
to satisfy the length condition $\ell(w_\ell)=\ell(w)+\ell(w_\ell^M)$, since long elements are all of order 2 (cf.~Section 5.1 of [CST]).

Next, let $w$ and $w'\in B(G)$.
Let $M_w$ and $M_{w'}$ be the associated Levi subgroups and $A_w$ and $A_{w'}$ their centers.
Set $A_w^{w'}=A_{w'}\cap M^d_{w'}$, where $M^d_{w'}$ is the derived group of $M_{w'}$.

Here we only need the case that $w=w_\ell$.
The product $A_{w_\ell}^{w'} A_{w'}\subset A_{w_\ell}=A$ is open of finite index in $A$ with $A_{w_\ell}^{w'}\cap A_{w'}$ a finite set.
Given $a$ in this product we have $a=bc$, $b\in A^{w'}_{w_\ell}$ and $c\in A_{w'}$.
Any other decomposition will be of the form $a=(b\xi^{-1})(\xi c)$ with $\xi\in A^{w'}_{w_\ell}\cap A_{w'}$, a finite set.
Moreover, write $c=c'z$ with $c'\in A'_{w'}=(A_{w'})'$ and $z\in Z$.
Here $A'_{w'}=(A_{w'})'$ is the splitting of $Z\backslash A_{w'}$ as in (6.3.1).

The following proposition, Proposition 5.7 of [CST], is the general version of Proposition 6.5.5.
Its proof takes a large part of [CST] and is inspired by the work of Jacquet and Ye [J1,JY], although our asymptotic expansion does not follow from their germ expansion.
In fact, no germ expansion is expected to exist for a partial Bessel function, and our asymptotics seems to be the best one can get.

\medskip\noindent
(6.5.7)\ {\bf Proposition}.
{\it Let $f\in\sM(\pi)$ with $W^f(e)=1$.
Then

(a)\ there exists $f_{1,e}\in C_c^\infty (G,\omega_\pi)$,

(b)\ for each $w'\in B(G)$, $w'\neq e$, there exists $f_{1,w'}\in C_c^\infty (\Omega_{w'}, \omega_\pi)$ such that for $\varphi$ with sufficiently large support depending on $f$, we have

(i)\ $B^G_\varphi (\dot w_\ell a,f)=B^G_\varphi (\dot w_\ell a, f_{1,e}) +\sum_{w'\in B(G)\atop w'\neq e} B^G_\varphi (\dot w_\ell a, f_{1,w'})$

\noindent
for all $a\in A$,

(ii)\ $B^G_\varphi (\dot w_\ell a, f_{1,e})$ depends on $\pi$ only through $\omega_\pi$ for all $a\in A$,

(iii)\ for each $w'\in B(G)$, $w'\neq e$, we have
$$
B^G_\varphi (\dot w_\ell a, f_{1,w'})=\omega_\pi (z) B^G_\varphi (\dot w_\ell b c', f_{1,w'})
$$
is uniformly smooth as a function of $c'\in A'_{w'}$, where $a=bc' z$.
}

\medskip\noindent
(6.5.8) {\bf Remark}.
Proof of this proposition is long and delicate.
One needs to show that the Bruhat cells for which $w\not\in B(G)$ do not contribute to the asymptotic expansion, and even for $w\in B(G)$, only its relevant part $C_r (\dot w)\colon = U\dot w A_w U$ contributes.
This is what happens for a full Bessel function as proved in [J1,JY].
In [CST] we show that the same can happen if we enlarge the support of defining $\varphi$ in $B^G_\varphi (-,f)$, depending on $f$.
This is the fist step of our argument in establishing the asymptotics of $B^G_\varphi (-,f)$, and is inspired by and follows the ideas in [J1,JY].
We refer the reader to the complete proof given in [CST] for next steps and details.

The arguments given for this asymptotics rely very heavily on the fact that $\sM(\pi)\subset C_c^\infty (G,\omega_\pi)$.
This is exactly why our proof of Theorem 3.5 was a reduction to a proof of (SCS) through our local--global deformation arguments.
In fact, a general direct proof of stability for any irreducible generic representation seems to be fairly out of reach since $\sM(\pi)$ is no longer contained in $C_c^\infty (G,\omega_\pi)$ and our arguments in [CST] cannot be applied.

\bigskip\noindent
7.\ {\bf Cases of exterior cube for $GL_n$.}
Among the cases appearing in the Langlands--Shahidi method is when $H=E_n^{sc}$ and $M_H^{der}=SL_n$, $n=6,7,8$, which gives the $L$--functions attached to the exterior cube representation of $GL_n(\bC)$, $n=6,7,8$.
The symbol $E_n^{sc}$ denotes the simply connected form of exceptional group $E_n$, $n=6,7,8$.

The case $n=6$ seems to be manageable and if $p\neq 2,3$, $p=\char(F)$, and already proves the equality of $L$--functions
$$
L(s,\Lambda^3 \rho)=L(s,\pi(\rho),\Lambda^3)\tag7.1
$$
for any continuous Frob--semisimple $n$--dimensional representation of $W'_F$.
The  equality of $\var$--factors and root numbers needs a proof of (SCS) in this case, whose analysis seems rather similar to $r=\Lambda^2$.
This is being studied by Cogdell, Shahidi and Varma.

When $n=7$ and $p\neq 7$ or $n=8$ and $p\neq 2$, the equality (7.1) is again valid.

The restrictions on $p$ in each case forces all the supercuspidal representations of $GL_n(F)$ to be monomial, $n=6,7,8$, and there is no need to consider virtual sums of monomial representations for $GL_n(F)$ with $n>8$ for which the equality (7.1) is not known.

The multiplicativity in this case comes from
$$
\Lambda^3 (\rho_1\oplus\rho_2)=\Lambda^3\rho_1\oplus\rho_1\otimes\Lambda^2\rho_2\oplus\Lambda^2\rho_1\otimes\rho_2\oplus\Lambda^3\rho_2.\tag7.2
$$
The equality of factors for $\rho_1\otimes\Lambda^2\rho_2$ and $\Lambda^2\rho_1\otimes\rho_2$, which are needed for proving (7.1), are also available through the cases appearing in the Langlands--Shahidi method.
In fact, for $n=6$, the only delicate case is when $R=\std\otimes\Lambda^2$ of $GL_2(\bC)\times GL_4(\bC)$.
The factors then come from the case $(D_{5,3})$ of [Sh8], pg.~188, in which we can take $H=GSpin(10)$ with $M_H=GL(2)\times GSpin(6)$ and identify $GSpin(6)$ with $GL(4)$ through isogeny.
For $n=7$, the case $R=\std\otimes\Lambda^2$ of $GL_m(\bC)\times GL_{7-m} (\bC)$, $m=2,3$, follows from the case $(E_{6,2})$ of [Sh8], pg.~188, when $m=2$, while for $m=3$, one can use the case $(D_{6,3})$ there.
Finally, for $n=8$, the case $R=\std\otimes\Lambda^2$ of $GL_m (\bC)\times GL_{8-m} (\bC)$, $m=2,3,4$, comes from the case $(E_{7,4})$ of [Sh8], pg.~190, when $m=2$, while the cases $m=3$ and 4 follow from cases $(E_{7,2})$ and $(D_{7,3})$ of [Sh8], respectively.
In the next section, we will address the equality of the factors in the case of $GL_n$ and for an arbitrary representation $r$, including a general discussion of multiplicativity.

\bigskip\noindent
8.\ {\bf The general case}.
We now address the general case for $GL_n$.
Let $\rho$ be as before an $n$--dimensional continuous Frobenius--semisimple representation of $W'_F$, where $F$ is a local field.
We now let $r$ be an arbitrary finite dimensional complex representation of $GL_n(\bC)$.
We then have the Artin factors $\var(s,r\cdot\rho,\psi)$, $L(s,r\cdot\rho)$ and $\gamma(s,r\cdot\rho,\psi)$, where
$$
\gamma (s,r\cdot\rho,\psi)=\var(s,r\cdot\rho,\psi) L(1-s,r^\vee\!\!\cdot\!\rho)/L(s,r\cdot\rho).\tag8.1
$$

Next let $\pi(\rho)$ be the irreducible admissible representation of $GL_n(F)$ attached to $\rho$ by LLC.
We recall from earlier sections that $\gamma(s,r\cdot\rho,\psi)$ satisfies stability in general and in particular (SCS), (AM) and (FE).
We will discuss multiplicativity/additivity in general in some detail soon.
But we fist assume the existence of a set of local factors $\var(s,\pi(\rho),r,\psi)$, $L(s,\pi(\rho),r)$ and $\gamma(s,\pi(\rho),r,\psi)$ satisfying (AM), (FE), (SCS) and (M), in which case further elaboration is needed as we explain later.
The question that we aim to discuss is 

\medskip\noindent
(8.2)\ {\bf Question}.
{\it Is the LLC robust enough to satisfy
$$
\gamma(s,r\cdot\rho,\psi)=\gamma(s,\pi(\rho),r,\psi)
$$
for all $r$?
Similarly for $\var$ and $L$.}

We will show that under the validity of (AM), (FE), (SCS) and (M), as to be refined and elaborated, all introduced in Section 4, the answer to question (8.2) is yes for every $r$.

Given $r$, let us call the existence of factors $L(s,\pi,r)$, $\var(s,\pi,r,\psi)$ and $\gamma(s,\pi,r,\psi)$ satisfying axioms (AM), (FE), (SCS) and (M), an $r$--theory.
(In [Sh9], this was called an $r'$--theory, with $r$--theory reserved for (AM), (FE), (M) and general stability axiom (S) instead of (SCS) in [Sh9].)

\medskip\noindent
(8.3)\ {\bf Multiplicativity/Additivity}.
To address the general form of multiplicativity, which is needed for an arbitrary representation $r$ of $GL_n(\bC)$, we need to recall a few facts and notions from finite dimensional representation theory of $GL_n(\bC)$.

\medskip\noindent
(8.3.1)\ {\bf Young diagrams and Schur functors}.
Let $m$ be a fixed positive integer.
Let $\lambda=(\lambda_1,\ldots,\lambda_m), \lambda_1\geq\lambda_2\geq\ldots\geq\lambda_m\geq 0,\ \lambda_i\in\bN \cup \{0\}$, be a partition of $\lambda_1+\ldots+\lambda_m$.
If $n\geq m$, we extend $\lambda$ to a partition of $n$ by letting $\lambda_{m+1}=\ldots=\lambda_n=0$.
Let $L_1,\ldots,L_n$ be the standard basis for $\bC^n$.
Then each partition $\lambda$ defines a unique finite dimensional irreducible representation of $SL_n(\bC)$ of highest weight 
$$
\lambda=\lambda_1 L_1+\ldots+\lambda_n L_n,\tag8.3.2
$$
which we denote by $\bS_\lambda(\bC^n)$.
If we set $\lambda_n=0$, then this will be a one--one correspondence between partitions and finite dimensional irreducible representations of $SL_n(\bC)$.
In fact, $\bS_\lambda(\bC^n)=\bS_{\lambda'} (\bC^n)$ if and only if $\lambda-\lambda'=(\alpha,\alpha,\ldots,\alpha),\alpha\in\bZ$.
We refer to [FH] for details.
The partition $\lambda=(\lambda_1,\ldots,\lambda_n)$ defines a Young diagram
$$
\vbox{
\offinterlineskip\tabskip=3pt
\halign{
\strut # &
\hfil # \hfil &
\vrule # &
\hfil # \hfil &
\vrule # &
\hfil # \hfil &
\vrule # &
\hfil # \hfil &
\vrule # &
\hfil # \hfil &
\vrule # &
\hfil # \hfil &
\vrule # &
\hfil # \hfil &
\vrule # &
\hfil # \hfil &
\vrule #\cr
\omit& &\multispan{15}{\hrulefill}\cr
& $\lambda_1$ & &  & &  & &  & &  & &  & &  & &  &\cr
\omit& &\multispan{15}{\hrulefill}\cr
& $\lambda_2$ & &  & &  & &  & &  & &  & &  &\cr
\omit& &\multispan{13}{\hrulefill}\cr
& $\cdot$ & &  & & & & & & & &\cr
\omit& & \multispan{9}{\hrulefill}\cr
& $\cdot$ & &  & & & & & &\cr
\omit& & \multispan{7}{\hrulefill}\cr
& $\cdot$ & &  & & & &\cr
\omit& & \multispan{5}{\hrulefill}\cr
}}
$$

\vskip-1truein
\leftline{(8.3.3)}

\vskip.75truein
\noindent
with $\lambda_i$ boxes in the $i$th row as in the picture in which case $\lambda=(7,6,4,3,2)$.
The symbol $\bS_\lambda$ is called the {\it Schur} functor.
It acts on $\bC^n$ for any $n$ which can be identified with the standard representation of $SL_n(\bC)$ on $\bC^n$, and thus takes the standard representation of $SL_n(\bC)$ and generates all its finite dimensional irreducible representations for any fixed $n$ as $\lambda$ varies.
One can check that $\bS_\lambda$ is indeed a functor on finite dimensional complex spaces $\bC^n$ as $n$ varies.

Now, if we define $a_i=\lambda_i-\lambda_{i+1},i=1,\ldots,n-1$, and assume $\lambda_n=0$ per our earlier comments, we get 
$$
\lambda=a_1 L_1+a_2 (L_1+L_2)+\ldots+a_{n-1} (L_1+\ldots+L_{n-1})\tag8.3.4
$$
with $L_1+\ldots+L_k$ the highest weight for the $k$th fundamental representation $\Lambda^k\bC^n$ of $SL_n(\bC)$.
Moreover, as in [FH], one sets
$$
\Gamma_{a_1,\ldots,a_{n-1}}\colon =\bS_\lambda (\bC^n)\tag8.3.5
$$
which appears in
$$
\Gamma_{a_1,\ldots,a_{n-1}}\subset \Sym^{a_1} (\bC^n)\otimes\ldots\otimes\Sym^{a_k}(\Lambda^k\bC^n)\otimes\ldots\otimes\Sym^{a_{n-1}} (\Lambda^{n-1}\bC^n)\tag8.3.6
$$

Finally, if $m>n$, then $L_{n+1}=\ldots=L_m=0$ and we only have $\lambda_1 L_1+\ldots +\lambda_n L_n$ and thus $\lambda_{n+1},\ldots,\lambda_m$ are irrelevant.

\medskip\noindent
(8.3.7)\ {\bf Young symmetrizer and Schur functor}.
Let $\lambda$ be a partition of length $|\lambda|=m$.
Let $S_m$ be the symmetric group in $m$ letters.
Define 
$$
P=P_\lambda=\{g\in S_m|g\text{ preserves each row of }Y_\lambda\},\tag8.3.8
$$
where $Y_\lambda$ denotes the Young diagram of $\lambda$, and 
$$
Q=Q_\lambda=\{g\in S_m|g\text{ preserves each column of }Y_\lambda\}.\tag8.3.9
$$
For example, if $m=3$, $\lambda=(2,1)$ giving
\vskip-.25truein
\qquad\qquad\qquad\qquad\qquad\qquad\qquad\quad \ \ $\vbox{
\offinterlineskip\tabskip=3pt
\halign{
\strut # &
\vrule # &
\hfil # \hfil &
\vrule # &
\hfil # \hfil &
\vrule #\cr
\omit&\multispan{5}{\hrulefill}\cr
& &  & &  &\cr
\omit&\multispan{5}{\hrulefill}\cr
& &  &\cr
\omit&\multispan{3}{\hrulefill}\cr
}}$ 

\vskip-.65truein

\hskip3.45truein for $Y_\lambda$, then $P_\lambda=\{1,(12)\}$ and

\bigskip\noindent $Q_\lambda=\{1,(13)\}$, where $(ab)$ means the transposition on $\{a,b\}$.

Let
$$
a_\lambda=\sum_{g\in P} e_g\tag8.3.10
$$
and
$$
b_\lambda=\sum_{g\in Q} (sgn\ g) e_g.\tag8.3.11
$$
Here $\{e_g |g\in G\}$ is a basis of the underlying vector space for the group algebra $\bC S_m$ of $S_m$, in which the multiplication is given by
$$
e_g\cdot e_h=e_{gh}.\tag8.3.12
$$

We now define the {\it Young symmetrizer} 
$$
c_\lambda\colon = a_\lambda \cdot b_\lambda\in\bC S_m.\tag8.3.13
$$

If $V$ is a complex vector space, we can let $S_m$ act on $V^{\otimes m}=\overbrace{V\otimes\ldots\otimes V}$ by permuting the factors.

We now let $GL(V)$ act on the left on $V^{\otimes m}$ diagonally, i.e.,
$$
g\cdot (v_1\otimes\ldots\otimes v_m)=g v_1\otimes\ldots\otimes g v_m,\tag8.3.14
$$
while $S_m$ acts on the right permuting the factors as we discussed, commuting with $GL(V)$.
The group algebra $\bC S_m$ of $S_m$ now acts on $V^{\otimes m}$.
In particular, we can consider the action of $c_\lambda$ on $V^{\otimes m}$.
We set 
$$
\bS_\lambda (V)\colon = \text{ Im} (c_\lambda|V^{\otimes m})\tag8.3.15
$$
on which $GL_n(\bC)$ will act if we realize $V=\bC^n$ and $GL(V)=GL(\bC^n)=GL_n(\bC)$, where $n=\dim_{\bC} V$.
The Schur functor is then the functor $\bS_\lambda$ attaching to any finite dimensional complex vector space $V$, the representation of $GL_n(\bC)$, $n=\dim_{\bC}V$, on $\text{Im}(c_\lambda| V^{\otimes m})$.
It is an irreducible representation of $GL_n(\bC)$ whose restriction to $SL_n(\bC)$ is of highest weight $\lambda$.

Let us now consider a few examples:

\medskip\noindent
(8.3.16)\ Let $\lambda=(m)$ for which the Young diagram is
$$
Y_{(m)}\colon 
\vbox{
\offinterlineskip\tabskip=3pt
\halign{
\strut # &
\vrule # &
\hfil # \hfil & 
\vrule # &
\hfil # \hfil &  
\vrule # &
\hfil # \hfil & 
\vrule # &
\hfil # \hfil & 
\vrule # &
\hfil # \hfil & 
\vrule #\cr
\omit& \multispan{11}{\hrulefill}\cr
& &  & &  & &  & &  & &  &\cr
\omit& \multispan{11}{\hrulefill}\cr
}}
$$
In this case $P_\lambda=P_{(m)}=S_m$, while $Q_\lambda=Q_{(m)}$ is just the identity and $c_{(m)}=\sum\limits_{g\in S_m} \ e_g$.
The image of every complex space $V$ is then generated by
$$
c_{(m)} (v_1\otimes\ldots\otimes v_m)=\sum_{\sigma\in S_m} v_{\sigma(1)}\otimes\ldots\otimes v_{\sigma(m)},
$$
the so called $m$th symmetric power of $V$, denoted by $\Sym^m V$.
It is irreducible and its restriction to $SL_n(\bC)$ has highest weight $(m)$. 

\medskip\noindent
(8.3.17)\ Let $\lambda=\oversetbrace{m}\to {(1,\ldots,1)}$ and thus

\vskip.5truein
$
\hskip2truein Y_{(1,\ldots,1)}\colon
$

\vskip-.85truein
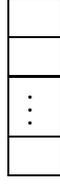
$$
\vbox{
\offinterlineskip\tabskip=2pt
\halign{
\strut # &
\vrule # &
\hfil # \hfil &
\vrule #\cr
\omit&\multispan{3}{\hrulefill}\cr
& &  &\cr
\omit&\multispan{3}{\hrulefill}\cr
& &  &\cr
\omit&\multispan{3}{\hrulefill}\cr
& & \ $\vdots$ &\cr
\omit&\multispan{3}{\hrulefill}\cr
& &  &\cr
\omit&\multispan{3}{\hrulefill}\cr
}}
$$

In this case $P_\lambda$ is identity, while
$$
b_\lambda=\sum_{g\in S_m}\text{ sgn} (g) e_g.
$$
The image of every complex space $V$ under $c_\lambda$ is then generated by
$$
\sum_{\sigma\in S_m}\text{ sgn}(\sigma) v_{\sigma(1)}\otimes\ldots\otimes v_{\sigma(m)},
$$
is the $m$th exterior power of $V$, denoted by $\Lambda^m V$.
It is irreducible and its restriction to $SL_n(\bC)$ has highest weight $(1,1,\ldots,1)$.
This is the $m$th fundamental representation of $SL_n(\bC)$.

\noindent
(8.3.18)\ Finally, let $m=3$ and $\lambda=(2,1)$ with $Y_\lambda\colon$
\vskip-.55truein
$$  
\qquad\qquad\qquad\qquad\ \ \ \qquad\vbox{
\offinterlineskip\tabskip=1pt
\halign{
\strut # &
\vrule # &
\hfil # \hfil &
\vrule # &
\hfil # \hfil &
\vrule #\cr
\omit&\multispan{5}{\hrulefill}\cr
& &  \ & & \ & \cr
\omit&\multispan{5}{\hrulefill}\cr
& &  \ &\cr
\omit&\multispan{3}{\hrulefill}\cr
}}$$

\vskip-.70truein
\hskip4.125truein{.}

\smallskip\noindent
In this case
$$
\aligned
c_{(2,1)}&= (e_1+e_{(12)})(e_1-e_{(13)})\\
&=1+e_{(12)}-e_{(13)}-e_{(132)}.
\endaligned\tag8.3.19
$$
Then the image of $c_{(2,1)}|V^{\otimes 3}$ is generated by
$$
v_1\otimes v_2\otimes v_3+v_2\otimes v_1\otimes v_3-v_3\otimes v_2\otimes v_1-v_3\otimes v_1\otimes v_2.\tag8.3.20
$$

\noindent
(8.3.21)\ {\bf Schur functors and multiplicativity}.
As explained before the problem of multiplicativity is that of how $\gamma$--factors behave under induction.
More precisely, we expect that for a finite dimensional representation $r$ of $GL_n(\bC)$ 
$$
\gamma(s,\pi_1\boxplus\pi_2,r,\psi)=\gamma(s,\pi_1,r,\psi)\gamma(s,\pi_2,r,\psi)\gamma(s,(\pi_1,\pi_2),R,\psi),\tag8.3.22
$$
whenever $\pi_i$ are unitary representations of $GL_{n_i}(F)$, $i=1,2$, and $\pi_1\boxplus\pi_2$ denotes the ``isobaric sum'' of $\pi_1$ and $\pi_2$ which is in fact a distinguished constituent of the representation of $GL_{n_1+n_2}(F)$ induced from the representation $\pi_1\otimes\pi_2$ on the Levi subgroup $GL_{n_1} (F)\times GL_{n_2}(F)$.
We will be more specific about ``isobaric sums'' in (8.3.34).
The representation $R$ of $GL_{n_1}(F)\times GL_{n_2}(F)$ is what we will determine using Schur functors as we have already seen when $r=\Lambda^2$ or $\Lambda^3$.

Let $\nu$ be a partition of length $|\nu|$ (partition of $|\nu|$) where $|\nu|$ is a positive number giving the length of $\nu=(\nu_1,\nu_2,\ldots)$, i.e., $|\nu|=\sum\limits_i\nu_i$.
Let $\bS_\nu$ be the Schur functor attached to $\nu$.
Then for any complex space $U$,
$$
\bS_\nu (U)=\text{ Im}(c_\nu |U^{\otimes |\nu|}),\tag8.3.23
$$
where $c_\nu$ is the corresponding Young symmetrizer discussed in (8.3.7).

Now, let $\lambda$ and $\mu$ be partititions with
$$
|\nu|=|\lambda|+|\mu|.\tag8.3.24
$$
Take positive integers $n_1, n_2$ and let $n=n_1+n_2$.
We then have Schur functors $\bS_\nu$, $\bS_\lambda$ and $\bS_\mu$ which we can apply to standard representation of $GL_n(\bC)$, $GL_{n_1}(\bC)$ and $GL_{n_2}(\bC)$, respectively.
We recall that as discussed in 8.3.1, $n$ and $|\nu|$ can be treated independent of each other.
Similarly for $n_1$ and $|\lambda|$ as well as $n_2$ and $|\mu|$.

Let $V$ and $W$ denote the standard representations of $GL_{n_1}(\bC)$ and $GL_{n_2}(\bC)$, respectively.
Then $GL_n(\bC)$ acts on $V\oplus W\simeq\bC^n$ and we have ([FH], p.~79--80):
$$
\bS_\nu (V\oplus W)=\bigoplus_{\lambda,\mu} N_{\lambda\mu\nu} (\bS_\lambda V\otimes\bS_\mu W).\tag8.3.25
$$
Here the sum runs over all the partitions of $|\lambda|$ and $|\mu|$ with $|\lambda|+|\mu|=|\nu|$ and the multiplicities $N_{\lambda\mu\nu}$ are non--negative integers given by the {\it Littlewood--Richardson rule}:\ \ The number of ways the Young diagram of $\lambda$ can be extended to the Young diagram of $\nu$ by a ``strict'' $\mu$--expansion (cf.~[FH], Appendix A, p.~455--456).
We note that 
$$
N_{\lambda 0\nu}=N_{\nu 0\nu} = N_{0\mu\nu}=N_{0\nu\nu}=1.\tag8.3.26
$$
Thus
$$
\bS_\nu (V\oplus W)=\bS_\nu V\oplus \bS_\nu W\oplus R_\nu (V,W),\tag8.3.27
$$
where $R_\nu (V,W)$ is a representation of $GL(V)\times GL(W)$.
We remark that $\bS_\nu V$ means the image of $V\subset V\oplus W$ under the functor $\bS_\nu$.

As an example, we note that when $r=\Lambda^2$, then
$$
R_\nu (V,W)=V\otimes W,\tag8.3.28
$$
while for $r=\Lambda^3$
$$
R_\nu (V,W)=\Lambda^2 V\otimes W\oplus V\otimes\Lambda^2 W.\tag8.3.29
$$

\noindent
(8.3.30) {\bf Remark}.
The Littlewood--Richardson rule is explained in [FH], together with the notion of strict $\mu$--expansion.
The basic case of this notion when $\mu=(m)$, Pieri formula, is explained in pages 79 (or 455) of [FH], equation (6.8) (or (A.7)).
The strict $\mu$--expansions of an arbitrary $\lambda$ are obtained by adding $m$ boxes to the Young diagram of $\lambda$, with no two in the same column.
For example, if $\nu=(4,2)$, $\mu=(4)$ and $\lambda=(1,1)$, then $\lambda$ cannot be expanded to $\nu$ by a strict $\mu$--expansion.
This basic case needs to be treated separately from the general case where $\mu\neq (m)$.

\medskip\noindent
(8.3.31)\ {\bf Example}.
We like to express $\bS_\nu(V\oplus W)$, where $\nu=(4,2)$, explicitly by means of (8.3.25).
We need to calculate $N_{\lambda\mu\nu}$.
We will give all partitions $\lambda$ and $\mu$, where $Y_\lambda$ can be expanded to $Y_\nu$ by a strict $\mu$--expansion.
We note that $|\lambda|+|\mu|=|\nu|=6$.
Possible pairs are $((4,2),0),\ ((2,1),(2,1)),\ ((2,1),(3))$, $((2,2),(2))$, $((3,2)),(1))$, $((3),(2,1))$, $((2),(2,2))$, $((1),(3,2))$ and $(0,(4,2))$.
The pairs $((1,1),(4))$ and $((4,(1,1))$ are not admissible since $Y_{(1,1)}$ cannot be expanded to $Y_{(4,2)}$ by a strict (4)--expansion and conversely.
Multiplicities of all pairs are 1.
For example, for the pair $((2,1),(2,1))$, although
\bigskip
\vskip-.65truein
$$\rightline{\vbox{
\offinterlineskip\tabskip=3pt
\halign{
\strut # &
\vrule # &
\hfil # \hfil &
\vrule # &
\hfil # \hfil &
\vrule # &
\hfil # \hfil &
\vrule # &
\hfil # \hfil &
\vrule #\cr
\omit&\multispan{9}{\hrulefill}\cr
& & $\ \ $ & & $\quad$ & & 1 & & 1 &\cr
\omit&\multispan{9}{\hrulefill}\cr
& & $\quad$ & & 2 &\cr
\omit&\multispan{5}{\hrulefill}\cr
}}\qquad\qquad\quad}$$
is a strict $(2,1)$--expansion of $(2,1)$ to $(4,2)$,
\vskip-0.5truein
$$
\hskip1.85truein\vbox{
\offinterlineskip\tabskip=3pt
\halign{
\strut # &
\vrule # &
\hfil # \hfil &
\vrule # &
\hfil # \hfil &
\vrule # &
\hfil # \hfil &
\vrule # &
\hfil # \hfil &
\vrule #\cr
\omit&\multispan{9}{\hrulefill}\cr
& & $\ \ \ $ & & $\quad$ & & 1 & & 2 &\cr
\omit&\multispan{9}{\hrulefill}\cr
& & $\quad$ & & 1 &\cr
\omit&\multispan{5}{\hrulefill}\cr
}}$$
is not since the list $2,1,1$ as prepared by going from right to left and top to bottom of the added boxes, does not satisfy the strict condition stated in the first paragraph of page 456 of [FH]. We thus have
$$
\gathered
\bS_{(4,2)}(V\oplus W)=\bS_{(4,2)}V\oplus(\bS_{(2,1)}V\otimes\bS_{(2,1)} W)\oplus\\ 
\noalign{\bigskip}
(\bS_{(2,1)} V\otimes\bS_{(3)} W)\oplus (\bS_{(2,2)} V\otimes\bS_{(2)} W)\oplus\\ 
\noalign{\bigskip}
(\bS_{(3,2)} V\otimes\bS_{(1)} W)\oplus (\bS_{(1)} V\otimes\bS_{(3,2)} W)\oplus\\ 
\noalign{\bigskip}
(\bS_{(2)} V\otimes \bS_{(2,2)} W)\oplus (\bS_{(3)} V\otimes\bS_{(2,1)} W)\oplus\bS_{(4,2)} W.
\endgathered
$$
We thus see that $R_{(4,2)}(V,W)$ is a direct sum of 7 tensor products and is fairly complicated compare to $\Lambda^2$ and $\Lambda^3$.

Let $\rho_i,i=1,\ldots,m$, be $m$ continuous Frobenius--semisimple $n_i$--dimensional complex representation of $W'_F$.
Let $\pi_i=\pi(\rho_i)$, $i=1,\ldots,m$, be the corresponding irreducible admissible representations of $GL_{n_i}(F)$, $n_i=\dim_{\bC}\rho_i$, $i=1,\ldots,m$.
Set $n=n_1+\ldots+n_m$.
Fix $N\in\bN$.
Let $r$ be the irreducible representation of $GL_N(\bC)$ defined by the partition $\nu$ of length $|\nu|$ and let $\bS_\nu (\bC^N)=r (GL (\bC^N))$.
As we discussed in Section 8.3.1, we may treat $|\nu|$ and $N$ as independent objects.
In particular, $\bS_\nu$ will apply on any $\bC^N$ the same way, depending only on $\nu$.

Using (8.3.27), we can write
$$
\bS_\nu (\bigoplus^m_{i=1} \rho_i)=\bigoplus^m_{i=1} \bS_\nu (\rho_i)\bigoplus R_\nu\cdot (\rho_1,\ldots,\rho_m),\tag8.3.32
$$
where $R_\nu$ is a representation of $GL_{n_1}(\bC)\times\ldots\times GL_{n_m}(\bC)$.

\medskip\noindent
(8.3.33) {\bf Arithmetic multiplicativity in general}.
For arithmetic (Artin) factors we have
$$
\gamma (s,r\cdot(\bigoplus^m_{i=1} \rho_i), \psi)=\prod^m_{i=1} \gamma (s,r\cdot\rho_i,\psi)\cdot \gamma(s,R\cdot (\rho_1,\rho_2,\ldots,\rho_m),\psi),
$$
which in view of the properties of arithmetic factors and (8.3.32) is a theorem.
Here $R=R_\nu$ with $r(GL(\bC^n))=\bS_\nu(\bC^n)$.

\medskip\noindent
(8.3.34) {\bf Analytic multiplicativity in general}.
Let $\pi_i, i=1,\ldots,m$, be irreducible admissible representations of $GL_{n_i}(F)$.
Analytic multiplicativity is more subtle.
First we need an irreducible admissible representation of $GL_n(F)$.
If $\pi_1\otimes\pi_2\otimes\ldots\otimes\pi_m$ is a quasi--tempered representation of the Levi subgroup $GL_{n_1}(F)\times\ldots\times GL_{n_m}(F)$ of $GL_n(F)$ with its complex parameter in the positive cone, then we can choose the unique Langlands quotient of Ind $\pi_1\otimes\ldots\otimes\pi_m$ as our choice.
On the other hand if $\pi_i$ are in addition $\psi$--generic with respect to a generic character of $U_{n_1}(F)\times\ldots\times U_{n_m}(F)$ defined by $\psi$, then we will choose the unique generic constituent of Ind $\pi_1\otimes\ldots\otimes\pi_m$ as our choice.
In both cases we will use $\boxplus^m_{i=1}\pi_i$ to denote this choice as we did in [Sh9].
We refer to Remark 8.3.47 for further discussion of this.

Analytic multiplicativity requires the existence of factors on $GL_{n_1} (F)\times\ldots\times GL_{n_m }(F)$ attached to $R=R_\nu$ such that ``analytic multiplicativity''
$$
\gamma (s,\boxplus^m_{i=1} \pi_i,r,\psi)=\prod^m_{i=1}\gamma(s,\pi_i,r,\psi)\gamma(s,(\pi_1,\ldots,\pi_m),R,\psi)\tag"{\it (M)}"
$$
holds.
The representation $R$ is a direct sum $\bigoplus\limits_j$ of tensor products of representations of $GL_{n_1}(\bC)\times\ldots\times GL_{n_m} (\bC)$ of the form $\bigotimes\limits_i \bS_{\nu_{ij}}(\bC^{n_i})$, where $\nu_{ij}$ are partitions of length strictly less than $|\nu|$, i.e., $|\nu_{ij}| < |\nu|$, for all $1\leq i\leq m$ and $j$, with possible multiplicities, but given completely explicitly by the general rules discussed concerning (8.3.27).

The factors $\gamma(s,(\pi_1,\ldots,\pi_m), R,\psi)=\gamma(s,\bigotimes\limits^m_{i=1} \pi_i,R,\psi)$ are expected to satisfy an $R$--theory from which one can deduce, along the same lines as for $r$ to be explained next, that 
$$
\gamma(s,R\cdot (\rho_1,\ldots,\rho_m),\psi)=\gamma(s,\bigotimes^m_{i=1}\pi_i,R,\psi),\tag8.3.35
$$
i.e., the answer to Question (8.2) for $R$ is positive.

\medskip\noindent
(8.3.36) {\bf Remark}.
In the generic case, i.e., within the Langlands--Shahidi method, given $r$, the corresponding $R$--theories always exist.
This is quite evident from all the examples we have seen so far (cf.~Section 7) and the general induction of the method.

We will now explain how the validity of an $r$--theory and its corresponding $R$--theories, and in fact only the stable version of (8.3.35) for each $R$, will lead to a proof of
$$
\gamma(s,r\cdot\rho,\psi)=\gamma(s,\pi(\rho),r,\psi)\tag8.3.37
$$
for any $\rho$, i.e., the answer to Question 8.2 is positive.

In fact, we will sketch how the arguments given in Section 5 for the case $r=\Lambda^2$ can be generalized to prove (8.3.37) for an arbitrary $r$.

We start with the following proposition which generalizes (5.1) to arbitrary $r$.

\medskip\noindent
(8.3.38)\ {\bf Proposition (stable equality)}.
{\it Fix $n\in\bN$ and let $\rho$ be an $n$--dimensional irreducible complex representation of $W_F$.
Let $r$ be as before a finite dimensional irreducible representation of $GL_n(\bC)$, given by a partition $\nu$ of length $|\nu|$.
Assume the validity of (FE), (AM) and (SCS) for $r$.
Moreover, assume the validity of (M) for all irreducible admissible representations $\pi_i$, $i=1,\ldots,m$, of $GL_{n_i}(F)$ with $\sum\limits^m_{i=1} n_i=n$ and for all partitions $(n_1,\ldots,n_i)$ of $n$, defining $\boxplus^m_{i=1} \pi_i$, i.e., the existence of factors $\gamma(s,(\pi_1,\ldots,\pi_m),R,\psi)$; as well as the validity of stable version of (8.3.35) for these factors, i.e., for $\gamma(s,\bigotimes\limits^m_{i=1} (\pi_i\otimes\chi), R,\psi)$, where $\chi$ is highly ramified.
Then
$$
\gamma(s,r\cdot(\rho\otimes\chi),\psi)=\gamma(s,\pi(\rho)\otimes\chi,r,\psi)\tag8.3.39
$$
for all highly ramified $\chi$.}

\medskip\noindent
{\bf Proof}.
As in Proposition 5.1, we first show the existence of a base point $(\rho_0,\pi(\rho_0))$ for which the stable equality holds for all $\chi\in\hat F^*$.
To proceed we appeal to Lemma 5.2.
Exactly as in the proof of Proposition 5.1, we choose a global representation $\tilde\rho$ of $W_k$ and a gr\"ossencharacter $\tilde\chi=\bigotimes\limits_w \tilde\chi_w$ of $k$ with $k_v=F$ and $\tilde\chi_v=\chi$, where $k$ is a number field.

We now use induction on $n$ and assume the validity of Proposition 8.3.38 for every local field $F$ and every $m<n$.
We first use $(M)$ and the stable version of (8.3.35), as well as our induction hypothesis, the stable equality (8.3.39) for all $m<n$, to conclude:
$$
\gamma(s,r\cdot(\tilde\rho_w\otimes\tilde\chi_w), \psi_w)=\gamma(s,\pi(\tilde\rho_w\otimes\tilde\chi_w), r,\psi_w)\tag8.3.40
$$
for all $w\neq v$, $w<\infty$.
We note that $\pi_i=\pi(\rho_i)$, $1\leq i\leq m$, when using $(M)$.
The equality (8.3.40) is also valid at all $w=\infty$ and for any $\chi_w$ by $(AM)$.
We now compare functional equations $(FE)$ for both $\tilde\rho\otimes\tilde\chi$ and $\pi(\tilde\rho\otimes\tilde\chi)\colon =\otimes_w \pi(\tilde\rho_w\otimes\tilde\chi_w)$, $\tilde\rho_w=\tilde\rho|k_w$, to conclude:

\medskip\noindent
(8.3.41)\ {\bf (Existence of a base point)}.
{\it There exists a pair $(\rho_0,\pi(\rho_0))$ with $\rho_0$ irreducible and thus $\pi_0$ supercuspidal, such that}
$$
\gamma(s,r\cdot(\rho_0\otimes\chi),\psi)=\gamma(s,\pi(\rho_0)\otimes\chi,r,\psi)\tag8.3.42
$$
for all $\chi\in\hat F^*$.

We now assume $\chi$ is highly ramified.
Then by arithmetic stability we have
$$
\gamma(s,r\cdot(\rho\otimes\chi),\psi)=\gamma(s,r\cdot(\rho_0\otimes\chi),\psi)\tag8.3.43
$$
for all highly ramified $\chi$, with ramification depending on $\rho$ and $\rho_0$.

We now appeal to (SCS) to conclude
$$
\gamma(s,\pi(\rho)\otimes\chi,r,\psi)=\gamma(s,\pi(\rho_0)\otimes\chi,r,\psi).\tag8.3.44
$$
Putting (8.3.42), (8.3.43) and (8.3.44) together we get
$$
\gamma(s,r\cdot(\rho\otimes\chi),\psi)=\gamma(s,\pi(\rho)\otimes\chi,r,\psi)\tag8.3.45
$$
for any irreducible $\rho$ with $\chi$ highly ramified, depending on $\rho$, completing Proposition 8.3.38.

\medskip\noindent
(8.3.46) {\bf Corollary}.
{\it Under the same assumptions Proposition 8.3.38 is valid for any $\rho$, i.e., not necessarily irreducible}.

\medskip\noindent
{\bf Proof}.
We need to use $(M)$.

\medskip\noindent
(8.3.47) {\bf Remark}.
We should point out that when it comes to $\gamma$--factors the choice of constituent of Ind $\pi_1\otimes\ldots\otimes\pi_m$ should be irrelevant and thus we may use $\boxplus^m_{i=1} \pi_i$ to denote this induced representation by itself.
We refer to multiplicativity for Rankin product $L$--functions in [JPSS], Theorem 3.1, pg.~404, as well as general multiplicativity in the context of Langlands--Shahidi method discussed in [Sh5,Sh8], as examples.

On the other hand this is not true if we consider $L$ and $\varepsilon$--factors for any individual constituent of the induced representation.

\medskip\noindent
(8.3.48) {\bf Proposition (equality for monomials)}.
{\it Assume $\rho$ is monomial, i.e., it is induced from a character of a subgroup of finite index in $W_F$.
Then under validity of (FE), (AM), (SCS) and (M), i.e., an $r$--theory, as well as validity of the stable version of (8.3.35), which can be deduced from the corresponding $R$--theory, we have}
$$
\gamma(s,r\cdot\rho,\psi)=\gamma(s,\pi(\rho),r,\psi).\tag8.3.49
$$

\medskip\noindent
{\bf Proof}.
Here we use the globalization of Harris [H], getting a $\tilde\rho$ and $\pi(\tilde\rho)$ as in the proof of Proposition 8.3.38 such that $k_v=F$ and $\tilde\rho|k_v=\rho$.
We then choose our gr\"ossencharacter $\tilde\chi=\otimes_w\tilde\chi_w$ such that $\tilde\chi_v=1$, but $\tilde\chi_w$ is highly ramified for every $w\neq v$, $w<\infty$, for which $\tilde\pi_w=\pi(\tilde\rho_w)$ is ramified.
We again  compare functional equations for $\pi(\tilde\rho)$ and $\tilde\rho$ exactly as in Proposition 8.3.38, using equation (8.3.40) at all places $w\neq v$, finite or infinite, using Corollary (8.3.46) and (AM), to conclude the proposition.

Having equality (8.3.49) of arithmetic and analytic factors for monomial representation, i.e., a basis for the Grothendieck ring of $W_F$, we now appeal to Brauer's theorem and (M) through (8.3.32), to conclude

\medskip\noindent
(8.3.50) {\bf Theorem}.
{\it Let $r$ be an irreducible representation of $GL_n(\bC)\simeq GL(\bC^n)$ given by a fixed partition $\nu$ as $r(GL_n(\bC))=\bS_\nu(\bC^n)$, satisfying an $r$--theory, as well as the stable version of (8.3.35) for corresponding representations $R$, both for all $n\in\bN$.
Then (8.3.49) is valid for all pairs $(\rho,\pi(\rho))$, i.e., the answer of Question 8.2 is positive.}

\medskip\noindent
(8.3.51) {\bf Remark}.
It is important to explain that to use Brauer's theorem we need to have the validity of $r$--theory and the stable version of (8.3.35) for {\bf all} $n$.
On the other hand if the residual character of $F$ does not divide $n$, then all the supercuspidal representations of $GL_n(F)$ will be monomial and Proposition (8.3.48) and $(M)$ will imply Theorem (8.3.50) and thus the equality of the factors for any $(\rho,\pi(\rho))$.
It is also important to treat $|\nu|$ and $n$ as independent numbers as we discussed earlier.
We will record this as follows.

\medskip\noindent
(8.3.52) {\it Assume $p\nmid n$, where $p=\char(\overline F)$.
Let $r$ be an irreducible representation of $GL_n(\bC)$, satisfying an $r$--theory, as well as the stable version of (8.3.35) for corresponding representations $R$.
Then (8.3.49) is valid for all pairs $(\rho,\pi(\rho))$.}

\medskip\noindent
(8.3.53) {\bf Remark}.
We refer to Section 7, where the cases of $r=\Lambda^3$ is discussed for $GL_n$, $n=6,7,8$, using Langlands--Shahidi method, in which case one gets complete results when $p\neq 2,3\  (GL_6)$, $p\neq 7\ (GL_7)$ and $p \neq 2\ (GL_8)$, i.e., in tame cases, without needing $\Lambda^3$ for higher $n>8$.

\bigskip\noindent
9.\ {\bf Braverman--Kazhdan/Ngo Program}.
In the last part of this paper we touch upon the theory of $\gamma$--factors as developed by Braverman and Kazhdan [BK], but in connection with Ngo's discussion [N1,N2] of Vinberg's theory of universal monoids [V].

Braverman--Kazhdan's theory of $\gamma$--factors is a generalization of the work of Godement and Jacquet for $GL_n$ (cf.~[GJ]), which we now briefly recall.

Let $F$ be a $p$--adic field and $\sigma$ an irreducible admissible representation of $GL_n(F)$.
In [GJ], Godement and Jacquet developed a theory of standard $L$--functions for $GL_n(F)$.
More precisely, these are the $L$--functions $L(s,\sigma,r)$ with $r$ the $n$--dimensional standard representation of the $L$--group $GL_n(\bC)$ of $GL_n$, attached to $\sigma$, generalizing the work of Tate [T1] on $GL_1$.
In particular, they developed a theory of $\gamma$--factors $\gamma(s,\sigma,\stan,\psi)$ satisfying the four axioms (M), (SCS), (FE) and (AM) we discussed in the first part.
Let us briefly recall how these $\gamma$--factors are defined.

These $L$--functions are obtained by means of certain zeta functions $Z(\Phi,s,f)$ whose terms we now explain.
The function $f$ is a matrix coefficient of $\sigma$ defined by a pair of vectors $v\in\sH(\sigma)$ and $\tilde v\in\sH(\tilde\sigma)$, where $\tilde\sigma$ is the contragredient of $\sigma$.
More precisely
$$
f(g)=\langle\sigma(g)v,\tilde v\rangle.
$$
We then define
$$
\overset\check\to f(g)=f(g^{-1})=\langle v,\tilde\sigma (g)\tilde v\rangle.
$$
Let $\sS(M_n(F))$ denote the space of smooth (locally constant) functions of compact support on $M_n(F)$, the space of $n\times n$ matrices with entries in $F$.
The second ingredient in the definition of the zeta function is a function $\Phi\in\sS (M_n(F))$.
The zeta function is then defined as
$$
Z(\Phi,s,f)=\int_{GL_n(F)} \Phi (x) f(x)|\det (x)|^s d^* x,
$$
where $s\in\bC$.

Next we recall the Fourier transform $\hat\Phi$ of $\Phi$ by
$$
\hat\Phi(x)=\int_{M_n(F)} \Phi (y)\psi(\trace (xy)) dy,
$$
where $1\neq \psi\in\hat F$ as before.
The measure $dy$ is normalized so that
$$
\Phi(0)=\int_{M_n(F)}\hat\Phi (y) dy.
$$
The $\gamma$--factors attached to $\pi$ and the standard representation $\gamma(s)=\gamma(s,\sigma,\stan,\psi)$ is defined by
$$
Z(\hat\Phi,n-s,\overset\check\to f)=\gamma(s)Z(\Phi,s,f)
$$
for all $f$ and $\Phi$.
It only depends on $s,\sigma$ and $\psi$.

In these notes we discuss only one aspect of Braveman--Kazhdan theory [BK], namely the generalization of the vector space $M_n(F)$ by means of Vinberg universal monoids [V].
We follow Ngo [N1,N2].
We do this as to relate [N1,N2] to [Sh8] in certain instances when $\gamma$--functions can also be obtained from Langlands--Shahidi method [Sh5,Sh8] and in particular when $r=\Lambda^2$ or $\Sym^2$ as in previous sections.

\noindent
{\bf (9.1)\ Monoids}.
We refer to the first few pages of [V] for the definition of a monoid.
Briefly, we start with an affine algebraic variety $S$ which has an associative multiplication which is a morphism
$$
\mu\colon S\times S\longrightarrow S\tag9.1.1
$$
of algebraic varieties, i.e., $S$ is a semigroup.
We will assume that $S$ is defined over $k$, an algebraically closed field which we further assume for simplicity is of characteristic zero.
$S$ can then be realized as a sub--semigroup of End$(V)$ for some vector space $V$ over $k$ (cf.~[V]).
If $S$ happens to have an identity, then the unit may be considered as the identity of $GL(V)$.

We will further assume that $S$ is irreducible as a variety over $k$.

\medskip\noindent
{\bf (9.1.2) Definition}.
A semigroup with an identity is called a {\it monoid}.
If $M$ is a monoid, then we use $G(M)$ to denote its subgroup of units.
If $M$ and $M'$ are two monoids and $\varphi\colon M\longrightarrow M'$ is a dominant morphism of algebraic semigroups, then $\varphi(G(M))$ is an open subgroup if $G(M')$ and thus $\varphi(G(M))=G(M')$.

An algebraic monoid $M$ is called {\it reductive} if $G(M)$ is a reductive group.
$G(M)$ is never semisimple unless $M$ is a group.

\medskip\noindent
{\bf (9.1.3) Example}.
Let $V$ be a vector space over $k$ and let $G'\subset GL(V)$ be a semisimple group.
Then 
$$
M=k^*\overline{G'}\subset\End (V)
$$
is a monoid with a one--dimensional center.
Such examples are of particular interest to [N1,N2] and us.

We will say a monoid is {\it normal} if it is normal as a variety, i.e., all its local rings are normal, meaning they are integrally closed domains.

Given a reductive monoid $M$, let $G=G(M)$ denote the group of its units.
Then $G\times G$ acts on $M$ by
$$
(g_1,g_2)\cdot m=g_1 m g_2^{-1}.
$$
Let $G'=G_{\der}$ be the derived group of $G$.
Set
$$
A\colon = M//(G'\times G')\tag9.1.4
$$
which we call the {\it abelianization} of $M$, for the invariant--theoretic quotient of $M$ under double action of $G'$ (cf.~[V]).

We will now further assume that the map
$$
\pi\colon M\longrightarrow A\tag9.1.5
$$
is flat (cf.~[V]).
Then the fibers of $\pi$ are equidimensional which we will assume are reduced.

If $T'$ is a maximal torus of $G'$, we let $T^+\colon=T'/Z'$, where $Z'$ is the center of $G'$.
Moreover, if $\{\alpha_1,\ldots,\alpha_r\}$ is the set of simple roots of $(T',B')$, a Borel pair, then $T^+=T^{\ad}$, the maximal torus of $G'/Z'$, can be identified with $\bG_m^r$ through the well--defined map
$$
t\mapsto (\alpha_1 (t),\ldots,\alpha_r (t)),\tag9.1.6
$$
$t\in T^{\ad}$.
We note that $r$ is the (semisimple) rank of $G'$.
Finally, set $$G^+=(T'\times G')/Z'.$$

From now on we will assume $G'$ is simply connected.
Let $\{\omega_1,\ldots,\omega_r\}$ be its set of fundamental weights.
They are simply defined by
$$
\kappa(\omega_i,\alpha^\vee_j)=\delta_{ij},
$$
where
$$
\alpha^\vee_j=2\alpha_j/\kappa (\alpha_j,\alpha_j)\tag9.1.7
$$
is the corresponding coroot and $\kappa$ is the Killing form.
Moreover, $\delta_{ij}$ is the Kronecker delta function.

Now, let $\rho_i$ be the fundamental weight of $G'$ attached to $\omega_i$ on the space $V_i$.
We will extend $\rho_i$ to $\rho_i^+\colon G^+\longrightarrow GL(V_i)$ by
$$
\rho_i^+ (t,g)=\omega_i (w_0 (t^{-1}))\rho_i (g),\tag9.1.8
$$
where $w_0$ is the long element of $W(G,T)\simeq W(G',T')$.
Here $t\in T'$ and $g\in G'$.
We now extend each $\alpha_i$ to $G^+$ by
$$
\alpha_i^+ (t,g)=\alpha_i (t)\tag9.1.9
$$
in which we may assume $t\in T^{\ad}$.
We thus get a map
$$
(\alpha^+,\rho^+)\colon G^+\longrightarrow \bG_m^r\times\prod^r_{i=1} GL(V_i)\tag9.1.10
$$
which is an embedding.

Vinberg's universal monoid $M^+$ (cf.~[N1,N2]) is the closure of $G^+$ in $\bA^r\times\prod^r_{i=1}\End (V_i)$ with $\bA=\bG_a$.
We note that $M^+$ depends only on $G'$ since so does $G^+$.

Let us now consider the exact sequence
$$
0 @>>> G' @>>> G @>>> T @>>> 0\tag9.1.11
$$
with $T$ the torus $G'\backslash G$, $G'=G_{\der}$.
Vinberg's universal monoids theory allows us to obtain a monoid $M$ with an open embedding $G\hookrightarrow M$ and with $G$ as the group of units of $M$.
We will then have
$$
\CD
&&&&0&& 0\\
&&&&@VVV @VVV\\
0@>>> G' @>>> G @>>> T @>>> 0\\
&&&&@VVV  @VVV\\
&&&&M @>\pi>> A @>>> 0,\\
\endCD\tag9.1.12
$$

\medskip\noindent
where $\pi$ is abelianization map discussed earlier.
Recall

$$
\CD
&&&&0&& 0\\
&&&&@VVV @VVV\\
0 @>>> G' @>>> G^+ @>>> T^+ @>>> 0,\\
&&&&@VVV  @VVV\\
&&&&M^+ @>\pi^+>> A^+ @>>> 0\\
\endCD\tag9.1.13
$$

\medskip\noindent
where $\pi^+$ is the abelianization map for $M^+$.

Vinberg's main theorem (Theorem 5 of [V]) simply states that
$$
M=A\times {}_{A^+}M^+\colon = \{(a,m^+)\in A\times M^+|\varphi_{ab} (a)=\pi^+ (m^+)\},\tag9.1.14
$$
the {\bf fibered product} of $A$ and $M^+$ over $A^+$, where
\smallskip
$$
\CD
T @>>> T^+\\
@VVV @VVV\\
A @>\varphi_{ab}>> A^+\\
@A\pi AA  @AA\pi^+A\\
M @>\varphi>> M^+,
\endCD\tag9.1.15
$$
\smallskip\noindent
giving $M$ as a closed subsemigroup of $A\times M^+$ with $\varphi$ the projection map on the second coordinate.
{\it In short, the monoid $M$ is completely determined by the map}
$$
\varphi_{ab}\colon A \longrightarrow A^+.\tag9.1.16
$$
Then
$$
G=T\times {}_{T^+}G^+.\tag9.1.17
$$

\noindent
{\bf (9.2)\ An important special case}.
We will now consider the case where $G'\backslash G\simeq\bG_m$, i.e., the torus $T\simeq G_m$.
We then have
\smallskip
$$
\CD
& & & & 0 & & 0\\
& & & & @VVV @VVV\\
0 @>>> G' @>>> G @>>> \bG_m @>>> 0,\\
& & & & @VVV @VVV\\
& & & & M @>\pi>> \bA^1
\endCD\tag9.2.1
$$
where $M$ is obtained from $M^+$ by $\varphi_{ab}$.

Before we proceed further we should point out that this is precisely the situation in which $G$ is a Levi subgroup of a maximal parabolic subgroup of a connected reductive group which is exactly the set up for the theory of $L$--functions developed via Langlands--Shahidi method.

Going back to our discussion we now have
$$
\CD
0 & & 0 &\\
@VVV @VVV\\
\bG_m @>>> T^+&=\bG_m^r\\
@VVV @VVV\\
\bA^1 @>\varphi_{ab} >> \bA^r
\endCD\tag9.2.2
$$
and $\varphi_{ab}$ restricts to $\lambda\colon\bG_m\to T^+$.
In fact, as discussed in [N2], every ``{\it dominant}'' cocharacter $\lambda\colon\bG_m\to T^+=T'/Z'=T^{ad}$ of $T^{ad}$, which simply means $\kappa(\lambda,H_{\alpha_i})\geq 0$ for every simple root $\alpha_i$ (see Appendix), extends to a morphism $\varphi_{ab}\colon\bA^1\longrightarrow\bA^r$, and thus determines a monoid $M$ having $G'$ as the derived group of corresponding group of invertible elements of $M$.

We will now explain how we can use this dominant cocharacter $\lambda$ of $T^{ad}$ to define a monoid $M^\lambda$ with the group of units $G^\lambda$ so that its complex dual $\hat G^\lambda$ will have a representation whose restriction to the derived group of $\hat G^\lambda$ has $\lambda$ as its highest weight.

We like to use [N2] and for that reason we use $G$ to denote the semisimple group $G'$ and we let $M^\lambda$ be the monoid attached to $\lambda$.
Then 
$$
\CD
M^\lambda & @>>> & M^+\\
@V\pi VV  & & @VVV\\
A=\bG_a & @>\varphi_{ab}>> & A^+=\bG_a^r
\endCD\tag9.2.3
$$
and $\lambda$ defines a homomorphism
$$
\theta_\lambda\colon \bG_m\longrightarrow \Aut (G)\tag9.2.4
$$
by
$$
\theta_\lambda (a)=\Int (\lambda (a)),\tag9.2.5
$$
where Int denotes conjugations by elements of $T^{ad}$ which is well--defined.
The group $G^\lambda$ of units of $M^\lambda$ is then
$$
\aligned
G^\lambda&=G\rtimes_\lambda\bG_m\\
&=\{(g,a)\in G\times\bG_m|(g,a)(g',a')\colon = (g\lambda (a) g'\lambda (a^{-1}), a a')\}
\endaligned\tag9.2.6
$$
In fact, restricting (9.2.3) to $G^\lambda$ we get 
$$
\CD
0 @>>> G @>>>G^\lambda @>\pi >>\bG_m @>>> 0\\
&&&&&& @VV\lambda V\\
&&&&&&\bG_m^r=T^{\ad}
\endCD\tag9.2.7
$$
which gives $G^\lambda$ as the semi--direct product of $G$ and $\bG_m$ through $\Int\cdot\lambda$ as it in fact splits our exact sequence.

As in [N2], let us now look at the dual setting.
The group $G$ being simply connected implies that $\hat G$ is adjoint.
Moreover $G^{\ad}=\hat G^{\sc}$.
The maximal torus $\hat T^{\sc}$ of $\hat G^{\sc}$ is dual to $T^{ad}$.
Let $Z$ be the center of $G$.
Then we have
$$
0 @>>> Z @>>> T @>>> T^{\ad} @>>> 0\tag9.2.8
$$
and
$$
0 @>>> \hat Z @>>> \hat T^{\sc} @>>> \hat T @>>> 0,\tag9.2.9
$$
where $\hat Z$ is the center of $\hat G^{sc}$.

The dominant cocharacter $\lambda$ of $T^{\ad}$ can be identified with a dominant character $\lambda\colon\hat T^{\sc}\to \bG_m$.
It will then be the highest weight of an irreducible representation
$$
\rho_\lambda\colon \hat G^{\sc} @>>> GL (V_\lambda).\tag9.2.10
$$
By irreducibility $\hat Z$ acts on $V_\lambda$ by scalars, inducing a morphism
$$
\omega_\lambda\colon \hat Z@>>> \bG_m.\tag9.2.11
$$
As we explain next we can use $\omega_\lambda$ to define a central extension
$$
0 @>>> \bG_m @>>> \hat G^\lambda @>>> \hat G @>>> 0\tag9.2.12
$$
together with a representation
$$
\rho_\lambda^+\colon \hat G^\lambda \longrightarrow GL(V_\lambda),\tag9.2.13
$$
where the central $\bG_m$ acts as homothety on $V_\lambda$.
We use two central exact sequences for $\hat G^{\sc}$ and $\hat G^\lambda$ and set up morphisms between them to define $\hat G^\lambda$.
Thus consider
$$
\CD
0 @>>> \hat Z @>>> \hat G^{\sc} @>>> \hat G @>>> 0\\
&& @VV\omega_\lambda V  @VV F V  @VV \text{id} V\\
0 @>>> \bG_m @>>> \hat G^\lambda @>>> \hat G @>>> 0,\endCD\tag9.2.14
$$
in which $\hat Z$ and $\bG_m$ are identified with their images in $\hat G^{\sc}$ and $\hat G^\lambda$, respectively.

The map $F$ is the covering map on $(\hat G^\lambda)_{der}$, the derived group of $\hat G^\lambda$, while $F|\hat Z=\omega_\lambda$.
Moreover, $\omega_\lambda(\hat Z)$ is the center of $(\hat G^\lambda)_{der}$ and thus $F$ is completely defined by $\omega_\lambda$.
We then have
$$
\hat G^\lambda=(\bG_m\times F(\hat G^{\sc}))/\omega_\lambda (\hat Z).\tag9.2.15
$$

The representation $\rho_\lambda^+$ is obtained by taking $F^{-1}(x)$, $x\in (\hat G^\lambda)_{\der}$ 
and applying $\rho_\lambda$ to it.
It is well--defined since the difference between two elements in the fiber $F^{-1}(x)$ is in the $\ker(\omega_\lambda)$.
Finally, $\bG_m$ acts by multiplication since $\omega_\lambda(\hat Z)$ does, and thus $\bG_m$ acts as homothety.
We note that $\hat G^\lambda=(G^\lambda)^\wedge$.

\medskip\noindent
{\bf (9.3)\ Examples and connections with Langlands--Shahidi method}.

First take $\hat G=PGL_n(\bC)$ and thus $G=SL_n$ and let $\lambda$ be the highest weight of the standard representation of $\hat G^{\sc}$.
Then $\ker(\omega_\lambda)=\{1\}$ and thus $\omega_\lambda(\hat Z)=\hat Z$.
Moreover, $F(\hat G^{\sc})=(\hat G^\lambda)_{\der}=SL_n$.
Thus by (9.2.15) $\hat G^\lambda=GL_n(\bC)$ and $G^\lambda=GL_n$.

Next again assume $\hat G=PGL_n(\bC)$ and thus $G=SL_n$.
Let $\lambda$ be the highest weights $\delta_2$ or $2\delta_1$ of exterior square $\Lambda^2$ or symmetric square $\Sym^2$ representations of $GL_n(\bC)$.
Note that $\hat Z=\langle\xi_n\rangle$, where $\xi_n$ is a primitive $n$--th root of 1.
Then the action of $\hat Z$ on $V_\lambda$ is $\xi_n\to\xi_n^2$.
If $n$ is odd, then $\ker(\omega_\lambda)=\{1\}$ since $\langle\xi_n\rangle=\langle\xi_n^2\rangle$.
Thus $G^\lambda=GL_n$.

Assume $n$ is even.
Then $\ker(\omega_\lambda)=\{1,\xi^{n/2}\}=\{\pm 1\}$.
Consequently
$$
(\hat G^\lambda)_{\der}=SL_n / \{\pm 1\}=\hat G_0,\tag9.3.1
$$
where $\hat G_0\colon =(\hat G^\lambda)_{\der}$.
This is the semisimple group
$$
SL_n\longrightarrow\hat G_0\longrightarrow P GL_n,\tag9.3.2
$$
whose character module is of index 2 in the weight lattice.
In this case
$$
\hat G^\lambda=(\bG_m\times\hat G_0)/\omega_\lambda(\hat Z).\tag9.3.3
$$
Note that $\omega_\lambda(\hat Z)$ is the center of $\hat G_0$, embedded in $\bG_m\times\hat G_0$ diagonally and through $\omega_\lambda$ in $\bG_m$.

One can then conclude (cf.~[Ki,Sh7]) that
$$
G^\lambda=(\bG_m\times SL_n)/S,\tag9.3.4
$$
where $S$ is the subgroup of squares of $Z(SL_n)=\langle\xi_n\rangle$, i.e., with $S=\langle\xi_n^2\rangle$.
We notice that this formulation works for both $n$ even or odd.
One then notices that $(G^\lambda)_{\der}=SL_n$, where $(\hat G^\lambda)_{\der}=SL_n/\{\pm 1\}$ when $n$ is even but $(\hat G^\lambda)_{\der}=SL_n$, otherwise.

For higher exterior and symmetric powers similar situation happens.
For example, for $\Lambda^m$ or $\Sym^m$, $\hat G_0$ will depend on what $d=(m,n)$ is.

Assume $m=3$.
If $3\nmid n$, then $\ker(\omega_\lambda)=\{1\}$ and $\hat G^\lambda=GL_n$.
If $3|n$, then $\ker(\omega_\lambda)=\langle\xi_n^{n/3}\rangle=\{1,\omega,\omega^2\}$, where $\omega^2+\omega+1=0$.
Thus $\hat G_0=SL_n/\langle\omega\rangle$ and
$$
\hat G^\lambda=(\bG_m\times\hat G_0)/\omega_\lambda(\hat Z)\tag9.3.5
$$
with $\omega_\lambda(\hat Z)=\langle\xi_n^3\rangle$.
We leave it to reader to determine what $G^\lambda$ is.

One can give a uniform explanation of $G^\lambda,\lambda=\delta_p$ or $\lambda=p\delta_i$ for $\Lambda^p$ and $\Sym^p$, respectively, when $p$ is a prime.
We will then always have
$$
G^\lambda=(\bG_m\times SL_n)/S,\tag9.3.6
$$
where $S$ is the subgroup of $p$--th powers in the center of $SL_n$.
Note that when $p\nmid n$, then $G^\lambda=GL_n$ for $\lambda=\delta_p$ and $p\delta_1$.

The group $G^\lambda$ in these cases are exactly the Levi subgroups of maximal parabolic subgroups of simply connected reductive groups $H^\lambda$ which within the Langlands--Shahidi method give the $L$--function $L(s,\pi,\lambda)$ for any irreducible admissible representation $\pi$ of $G^\lambda(k)$.
We refer to [L1] and [Sh4] for the cases of $\lambda=\delta_m+ (-\delta_n)$, when 
$$
G^\lambda=SL_{m+n}\cap (GL_m\times GL_n),\tag9.3.7
$$
which gives the Rankin product $L$--function $L(s,\sigma_1\times\tilde\sigma_2)$ for $GL_m(k)\times GL_n(k)$, as well as cases of $\lambda$ giving $\Lambda^2$ and $\Sym^2$ for $GL_n(\bC)$ and $\Lambda^3$ for $GL_n(\bC)$, $n=6,7,8$.

At the end of this section we will look at the cases of symmetric power $L$--functions for $GL_2$ and compute the corresponding monoids as well as their groups of units in detail.
They agree with our discussion above.

We will now go to the general setting of a complex Kac--Moody group and explain how the Levi subgroups of maximal parabolic subgroups are the same as the groups $\hat G^\lambda$ for any dominant $\lambda$.
Thus within the conjectural generalization of Langlands--Shahidi method all the groups $\hat G^\lambda$ appear as the Levi subgroups of certain Kac--Moody groups, since all the irreducible finite dimensional complex representations of a complex reductive group, satisfying our condition (9.2), appear as a subrepresentation of the adjoint action of the complex group as a Levi subgroup of a maximal parabolic subgroup, on the Lie algebra of the unipotent radical of that parabolic subgroup.

More precisely, the structure theory of Kac--Moody groups says that there exists a choice $(\tilde H,\tilde G)$ of a complex adjoint Kac--Moody group $\tilde H$ and a Levi subgroup $\tilde G$ of a maximal parabolic subgroup $\tilde P=\tilde G\tilde N$ such that
$$
\hat G^{\sc}\longrightarrow (\tilde G)_{\der} @> r >> GL(V_1)\tag9.3.8
$$
gives the representation $V_\lambda$ of $\hat G^{\sc}$ defined by the highest weight $\lambda$ as a subrepresentation.
Here $r$ is the adjoint action of $\tilde G$ on the Lie algebra Lie$(\tilde N)$ of $\tilde N$ and as it is standard in Langlands--Shahidi method [Sh5,Sh8], $r$ decomposes as
$$
r=\bigoplus_i r_i\tag9.3.9
$$
with the order given in [Sh5,Sh8] and where $V_i$ is the space of $r_i$.
The fact that every $V_\lambda$ can be obtained as a subrepresentation of $r_1$ is the main induction of Langlands--Shahidi method and its natural extension to Kac--Moody groups [K].

The center $\hat Z$ of $\hat G^{\sc}$ will map onto the center of $(\tilde G)_{\der}$ and by Lemma 4.8 of [Sh8] its action on $V_1$ is simply multiplication and thus as homothety.
It will then be according to character $\omega_\lambda$ of $\hat Z$ discussed earlier since $(r_1,V_1)$ contains $V_\lambda$ as a subrepresentation.
In particular, $\omega_\lambda (\hat Z)$ is exactly the center of $(\tilde G)_{\der}$.
By construction
$$
\tilde G=(\bG_m\times (\tilde G)_{\der})/\omega_\lambda (\hat Z)=\hat G^\lambda.\tag9.3.10
$$

We record this discussion as 

\medskip\noindent
{\bf (9.3.11)\ Proposition}.
{\it Let $\lambda$ be a dominant cocharacter of $T^{\ad}$.
Denote by $G^\lambda$ the group of units of monoid $M^\lambda$ attached to $\lambda$ by Vinberg's universal monoids theory.
Let $\tilde H$ be a complex adjoint Kac--Moody group and $\tilde G$ a Levi subgroup of a maximal parabolic subgroup $\tilde P=\tilde G\tilde N$ of $\tilde H$ such that the adjoint action $r$ of $\tilde G$ on Lie $(\tilde N)$ decomposes as $r=\bigoplus\limits_i r_i$ with $r_1\cdot\eta$ containing
$V_\lambda$ as a subrepresentation, where $\eta\colon\hat G^{\sc}\longrightarrow (\tilde G)_{\der}$ is the covering map and $V_\lambda$ is the representation of $\hat G^{\sc}$ with highest weight $\lambda$.
Then $\tilde G\simeq\hat G^\lambda$.}

We now state and prove the needed result from Kac--Moody theory.
The version we present here was provided to us by Steve Miller which he calls a ``Folklore'' as it may be known to others, and in particular to Braverman.

\medskip\noindent
{\bf (9.3.12)  Proposition}.
{\it  Let $\rho$ be an irreducible finite dimensional complex representation of a simply connected complex reductive Lie group $G$ for which $G/G_{\der}$ is one dimensional.
Then there exists a Kac--Moody group $H$ and a maximal parabolic subgroup $P\subset H$ with a Levi decomposition $P=LN$ with $L_{\der}=G_{\der}$ such that $\rho$ appears in the adjoint action of $L$ on Lie$(N)$ and more precisely in $r_1$.}

\medskip\noindent
{\bf Proof}.
Let $\{\alpha_i\}$ be a set of simple roots for $G$ and denote by $\{\omega_i\}$ the set of fundamental weights of $G$.
Let $\lambda$ be the highest weight of $\rho$.
Choose non--negative integers $v_i$ such that
$$
\lambda=\sum_i v_i \omega_i.\tag9.3.13
$$
Write $v=(v_1,v_2,\ldots)$.
Let $C_0$ be the Cartan matrix of $G_{\der}$ and denote by $C$ the generalized Cartan matrix
$$
C=\bmatrix C_0 & -^t v\\
-v&2\endbmatrix.\tag9.3.14
$$
Let $H$ be a Kac--Moody group defined by $C$ and let $P=P_\alpha=LN$ be the maximal parabolic subgroup corresponding to the submatrix $C_0$, where $\alpha$ is the unique simple root for which the root vector $X_\alpha$ sits in Lie$(N)$.

For each $i$, let $H_{\alpha_i}\in \frak t$ be the semisimple member of $sl_2$--triple $(H_{\alpha_i},X_{\alpha_i},Y_{\alpha_i})$.
Let $X_\alpha$ be the corresponding root vector in Lie$(N)$ attached to $\alpha$.
Here $\frak t=$ Lie$(T)$, where $T$ is a maximal torus in $L$ and thus $H$.
Then by equation (0.3.1) of [K]
$$
\aligned
\text{Ad}(H_{\alpha_i}) X_\alpha&=[H_{\alpha_i},X_\alpha]\\
&=\alpha(H_{\alpha_i}) X_\alpha\\
&=\langle\alpha_i^\vee,\alpha\rangle X_\alpha\\
&=c_{in} X_\alpha,
\endaligned\tag9.3.15
$$
where $n$ is the size of the matrix $C=(c_{ij})$.

Now consider (9.3.13) again, i.e.,
$$
\lambda=\sum v_i\omega_i.
$$
Then
$$
\langle\lambda,\alpha^\vee_\ell\rangle =\sum v_i \langle \omega_i,\alpha_\ell^\vee\rangle=v_\ell,\tag9.3.16
$$
where $v_\ell$ is from amongst $v_i$.
Thus
$$
-c_{in}=v_i=\langle\lambda,\alpha_i^\vee\rangle=\lambda (H_{\alpha_i}).\tag9.3.17
$$
Using (9.3.15) and (9.3.17) it is now clear that under the adjoint action $\frak t$ acts on $X_\alpha$ by the lowest weight character $-\lambda$.
In particular $\rho$ appears in $r_1$, where $r_1$ is as in graded decomposition (9.3.9) of the adjoint action $r$.

\medskip\noindent
{\bf (9.3.18)\ Symmetric powers for $GL_2$}.
In this paragraph we will compute the monoids $M^\lambda=M_n$ that are defined by $\lambda=n\delta_1$ and thus give $n$--th symmetric power $L$--functions for $GL_2$.
Here $\delta_1$ is the first (and only) fundamental representation of $SL_2(\bC)$.
We first note that as discussed before $G^\lambda$, $\lambda=n\delta_1$, can be easily calculated to be $G^\lambda=GL_1\times SL_2$ when $n$ is even, while $G^\lambda=GL_2$ when $n$ is odd.
This agrees completely with the cases $\Sym^2$ and $\Sym^3$ of $GL_2$ which show up for the pairs $(H^\lambda,G^\lambda)=(Sp_4,GL_1\times SL_2)$ for $\Sym^2$ and $(H^\lambda,G^\lambda)=(G_2,GL_2)$ for $\Sym^3$ in the lists [L1,Sh4], both in the simply connected setting.
But knowledge of $G^\lambda$ up to these isomorphism will not allow us to compute the attached monoids which we now calculate per our earlier discussion.
It is instructive to do it most formally.

We will use our earlier notation and thus $G'=SL_2$.
With notation as in (9.1.8)--(9.1.10), we have 
$$
\aligned
\rho_1^+ (t,g)&=\omega_1 (w_0 (t^{-1}))\rho_1 (g),\\
&=a\rho_1(g)
\endaligned
$$
where $\omega_1=\delta_1$ is our only fundamental weight and $t=\pmatrix a&0\\ 0&a^{-1}\endpmatrix\in SL_2$.
Moreover 
$$
\alpha_1^+ (t,g)=\alpha(t)=a^2.
$$
The embedding from $G^+$ into $\bG_m\times GL(\bA^2)$, given by
$$
(\alpha^+,\rho^+)\colon (t,g)\mapsto (a^2,ag)
$$
will then give $M^+=\End (\bA^2)$ as expected from the case of standard $L$--functions for $GL_2$ of Godement--Jacquet.

To compute monoid $M_n$ from definition (9.1.14), we need to determine what $\varphi_{ab}$ and $\pi^+$ are in the case $\lambda=n\delta_1$.
Recall that $n\delta_1$ was the highest weight of $\Sym^n$, a representation of $SL_2(\bC)$, and thus a character of $GL_1(\bC)$ as the maximal torus of $SL_2(\bC)$.
We will use the same notation $n\delta_1$ to give the cocharacter $n\delta_1\colon\bG_m\to\bG_m$.
It is simply $(n\delta_1)(x)=a^n$.

The cocharacter $n\delta_1$ being dominant will extend to a morphism which is our $\varphi_{ab}\colon A_n\to A^+$ of definitions (9.1.14) and (9.1.15) with $A_n$ and $A^+$ abelianizations of $M_n$ and $M^+$, respectively.

Next we need to calculate $\pi^+\colon M^+\to A^+$.
Recall that $A^+$ is the invariant--theoretic quotient of $M^+$ under $G'\times G'$ as in (9.1.4) and (9.1.5).
Then, if $Z^+$ is the center of $GL_2$ inside $M^+$, $\pi^+ (\overline{Z^+})=A^+$ by part 2) of Theorem 3 of [V], where $\overline{Z^+}$ is the closure of $Z^+$.

Note that the determinant map
$$
M^+ \overset\det\to\longrightarrow \bA^1
$$
factors through $\bA^+$ since kernel of $\pi^+|GL_2$ is exactly $G'=SL_2$ by part 1) of the same Theorem 3 of [V], and the induced map $A^+\simeq\bA^1$ is an isomorphism.
Thus $\pi^+=\det$.

Now by definition (9.1.14)
$$
\aligned
M_n&=\{(a,m^+)\in A_n\times M^+\ |\ (n\delta_1)(a)=\pi^+(m^+)\}\\
&=\{(a,m^+)\in\bA^1\times\End (\bA^2)\ |\ a^n=\det (m^+)\}.
\endaligned\tag9.3.19
$$
Then the group of units
$$
G_n=\{(a,g)\in\bG_m\times GL_2\ |\ a^n=\det g\}\tag9.3.20
$$

We now appeal to the following to conclude that as abstract groups $G_n\simeq GL_1\times SL_2$, if $n$ is even, and $G_n\simeq GL_2$, otherwise.

\medskip\noindent
{\bf (9.3.21)\ Lemma}.
{\it Let $n$ be a positive integer and define
$$
G_n=\{(a,g)\ |\ a^n=\det g\}\subset GL_1\times GL_2.
$$
Then $G_n=GL_1\times SL_2$ if $n$ is even, and $G_n=GL_2$, otherwise.}

\medskip\noindent
{\bf Proof}.
Assume $n$ is even.
Write $n=2\ell$.
Define the map
$$
GL_1\times SL_2\longrightarrow G_n
$$
by
$$
(a,g_1)\mapsto (a,a^\ell g_1),
$$
where $a^\ell\colon = \pmatrix a^\ell&0\\ 0&a^\ell\endpmatrix$.
Note that this map is an isomorphism.

Now assume $n$ is odd.
Write $n=2\ell+1$.
Then the map
$$
G_1\longrightarrow G_n
$$
defined by
$$
(a,g)\mapsto (a,a^\ell g)
$$
is an isomorphism for all odd $n$.
But note that by definition
$$
G_1=\{(a,g)|\det g=a\}\subset GL_1\times GL_2
$$
is just $GL_2$, completing the lemma.

\medskip\noindent
{\bf (9.4)\ A Fourier transform}.
One of the main ingredients in defining the local coefficients and thus $\gamma$--factors within the Langlands--Shahidi method is intertwining operators which we briefly recall.
We refer to Chapter 4 of [Sh8] for details and references.
We will change our notation from $H$ to $G$ to agree with standard references.

Let $G$ be a connected reductive group over a local field $F$.
Fix a minimal parabolic subgroup $P_0=M_0 N_0$ over $F$.
Let $P$ be a parabolic subgroup of $G$, containing $P_0$ and thus standard with respect to $P_0$.
Fix a Levi decomposition $P=MN$, with $N\subset N_0$ and $M\supset M_0$.
Let $A$ be a maximal split torus of $G$ in $M_0$ and denote by $W(G,A)$ the Weyl group of $A$.
Choose a $w\in W(G,A)$ such that $w$ sends simple roots of $M$ into simple roots.

Let $\sigma$ be an irreducible admissible representation of $M(F)$.
For simplicity we will build the usual complex parameter $\nu$ into the central character of $\sigma$, and make the convergence of the operator depending only on the central character of $\sigma$, thus avoid mentioning $\nu$ (cf.~[Sh8] for details on convergence).
The intertwining operator $A(\sigma,w)$ is formally defined by:
$$
A(\sigma,w) f(g)\colon =\int_{N_w(F)}  f(\dot w^{-1} n g) dn,\tag9.4.1
$$
in which the subgroup $N_w\subset N$ is
$$
N_w\colon = N_0\cap \dot w N^- \dot w^{-1}\tag9.4.2
$$
which can be also identified with
$$
N_w\simeq N'\cap \dot w N \dot w^{-1}\backslash N',\tag9.4.3
$$
where $N'\subset N_0$ is the unipotent radical of the standard parabolic for which $w(M)$ is a Levi subgroup.
Here $N^-$ is the opposite group to $N$ and the element $\dot w$ is a representative for $w\in W(G,A)$ in $G(F)$.
The function $f\in V(\sigma)$, the space of representation $I(\sigma)$ induced from $\sigma$.
We note that $A(\sigma,w)$ intertwines $I(\sigma)$ and $I(w(\sigma))$.

We will now assume $P$ is self--associate and thus its opposite parabolic $\overline P=M\overline N$, $\overline N=N^-$, is conjugate to $P$ by a $w_0\in W(G,A)$.
More precisely, we take $w_0=w_\ell\cdot w_{\ell,M}^{-1}$, where $w_\ell$ and $w_{\ell,M}$ are long elements of $W(G,A)$ and $W(M,A)$, respectively.
Then $w_0(M)=M$ and $N_{w_0}=N$.
We will first interpret this operator as a convolution operator on $\overline N(F)$.
Write
$$
\dot w_0^{-1} n\overline n=m n',\tag9.4.4
$$
$n, n'\in N(F)$, $m\in M(F)$, $\overline n\in\overline N(F)=N^- (F)$.
This is slightly different from the usual way presented in [Sh6]
and is valid for an open dense subset of $N(F)$.
Then
$$
\aligned
A(\sigma,\dot w_0)f(g)=&\int_{N(F)} f(\dot w_0^{-1} ng)dn\\
=&\int_{N(F)} \sigma(m) f((\overline n)^{-1} g) dn.
\endaligned\tag9.4.5
$$
If we now restrict to $\overline N$, then
$$
A(\sigma,\dot w_0) f(\overline n_1)=\int_{N(F)}\sigma(m) f((\overline n)^{-1}\overline n_1) dn\tag9.4.6
$$
which can be written as $(\Phi * f)(\overline n_1)$, the convolution of $f\in V(\sigma)$ with the measure 
$$
\Phi=\sigma(m) dn,\tag9.4.7
$$
i.e.,
$$
\aligned
(\Phi * f)(\overline n_1)=&\int\sigma(m) f((\overline n)^{-1} \overline n_1)\ dn\\
\noalign{\smallskip}
=&\Phi(\overline n\longrightarrow f((\overline n)^{-1} \overline n_1)),
\endaligned\tag9.4.8
$$
which is precisely how the convolution with a measure or distribution acts.

We now assume $G$ is quasisplit.
Then $P_0$ becomes a Borel subgroups $B$ of $G$ and $M_0=T$ a maximal torus of $G$ with $T\supset A$.
Moreover, we will assume $\sigma$ is generic and show that applying the canonical Whittaker functional of $V(w_0(\sigma))$, (cf.~[Sh8]), to $A(\sigma,w_0)$, behaves exactly as a Fourier transform does.

More precisely, let $\psi$ be a non--trivial additive character of $F$.
Then $\psi$ together with the splitting $(G,B,T,\{x_\alpha\}_\alpha)$ of our group (cf.~[Sh8]), defines a generic character of $U(F)$ which we still denote by $\psi$.
It will also give a set of representatives for every $w\in W(G,A)$ (cf.~Remark 8.2.1 of [Sh8]).
Now, assume $\sigma$ is $\psi$--generic.
Let $\lambda$ and $\lambda'$ be the canonical Whittaker functionals for $V(\sigma)$ and $V(w_0(\sigma))$, respectively.
The local coefficient $C_\psi(\sigma)$ is then defined by
$$
C_\psi(\sigma)^{-1}\lambda=\lambda'\cdot A(\sigma,\dot w_0)\tag9.4.9
$$
(cf.~[\quad,\quad]).

By definition in which $\psi(\on_1)\colon=\psi(w_0^{-1}\on_1 w_0),\on_1 \in \oN(F)$,
$$
\lambda'(A(\sigma,\dot w_0) f)=\int_{\overline N(F)} \lambda_M (A(\sigma,\dot w_0) f(\overline n_1)) \psi (\overline n_1) d\overline n_1,\tag9.4.10
$$
where $\lambda_M$ is a Whittaker functional for $V(\sigma)$.
Then, at least formally,
$$
\gathered
\lambda'(A(\sigma,\dot w_0) f)=\int_{\overline N(F)} \lambda_M ((\Phi * f)(\overline n_1)) \psi(\overline n_1)d\overline n_1\\
\noalign{\smallskip}
=\int_{\overline N(F)}\int_{N(F)} \lambda_M (\sigma(m) f((\overline n)^{-1} \overline n_1)) \psi(\overline n_1)dn d\overline n_1\\
\noalign{\smallskip}
=\int_{\overline N(F)} \lambda_M (( \int_{N(F)} \sigma(m)\psi(\overline n) dn) f (\overline n_1)) \psi (\overline n_1) d\overline n_1\\
\noalign{\smallskip}
=\int_{\overline N (F)} \lambda_M (\psi(\Phi) f (\overline n_1)) \psi (\overline n_1) d\overline n_1\\
\noalign{\smallskip}
=\lambda_M (\psi(\Phi) \int_{\overline N(F)} f(\overline n_1) \psi (\overline n_1) d\overline n_1)\\
\noalign{\smallskip}
=\lambda_M (C_\psi (\sigma)^{-1}\int_{\overline N(F)} f (\overline n_1) \psi (\overline n_1) d\overline n_1)
\endgathered\tag9.4.11
$$
by (9.4.9).

Here
$$
\psi(\Phi)\colon =\int_{N(F)} \sigma(m) \psi (\overline n) dn\tag9.4.12
$$
is the Fourier transform of the measure $\Phi$ defined by (9.4.7).

Since $f$ is of compact support modulo $P$, we can choose it so that
$$
u=\int_{\overline N(F)} f(\overline n_1) \psi(\overline n_1) d\overline n_1\tag9.4.13
$$
becomes any arbitrary vector in $V(\sigma)$.
It now follows from (9.4.11) that
$$
\lambda_M ((\psi(\Phi)-C_\psi(\sigma)^{-1}) u)=0
$$
for all $u\in V(\sigma)$, or 
$$
\psi(\Phi)\equiv C_\psi (\sigma)^{-1} \text{mod } (\ker (\lambda_M)).\tag9.4.14
$$
It can be easily checked that $\psi(\Phi)$ commutes with the action of $U_M(F)$.
We record this as:

\medskip\noindent
{\bf (9.4.15)\ Proposition}.
{\it Let $\psi(\Phi)$ be the $\psi$--Fourier transform of the measure $\sigma(m)dn$, i.e.,
$$
\psi(\Phi)=\int_{N(F)} \sigma(m)\psi(\overline n)dn.
$$
Then
$$
\psi(\Phi)\equiv C_\psi (\sigma)^{-1} \text{mod }(\ker (\lambda_M)).
$$
Moreover, the operator $\psi(\Phi)$ commutes with the action of $U_M(F)$, where $U_M=U\cap M$.}

\medskip\noindent
{\bf (9.4.16)\ Remark}.
Let $G$ be a quasisplit connected reductive group over a number field $k$ whose ring of adeles is $\bA_k$.
Let
$$
\rho\colon ^L\!G\longrightarrow GL_n (\bC)\times\Gamma_k
$$
be an analytic representation of $^LG$.
Let $\pi=\otimes_v\pi_v$ be a cuspidal automorphic representation of $G(\bA_k)$.
Let $L(s,\pi_v,\rho_v)$ be the local $L$--function attached to $\pi_v$ and $\rho_v$ by Langlands, whenever $\pi_v$ and $G$ as a group over $k_v$, are unramified.
We shall now assume that we have a theory of $L$--functions for $\rho$, i.e., a collection of $L$--functions $L(s,\pi_v,\rho_v)$ and root numbers $\var(s,\pi_v,\rho_v,\psi_v)$, $\psi=\otimes_v\psi_v$ a character of $k\backslash\bA$, $\psi\neq 1$, such that the global objects
$$
L(s,\pi,\rho)=\prod_v L(s,\pi_v,\rho_v)
$$
and
$$
\var(s,\pi,\rho)=\prod_v\var (s,\pi_v,\rho_v,\psi_v)
$$
have meromorphic continuation to all of $\bC$, satisfying the functional equation
$$
L(s,\pi,\rho)=\var(s,\pi,\rho) L(1-s,\tilde\pi,\rho).
$$

In [La], Lafforgue introduces the notion of a ``function of type $L$'' on $G(\bA_k)$ so as to relate the theory of $L$--functions to spectral analysis on $L^2(G(k)\backslash G(\bA_k))$.
The local $L$--functions and root numbers allow the definition of a Fourier transform for a function of type $L$ and the global functional equation for $L(s,\pi,\rho)$ then implies a non--linear Poisson summation formula for such functions.
He then shows how this non--linear Poisson summation formula can be used to build a kernel function to implement functoriality for $\rho$.
In particular, he concludes that such Poisson summation formulas are equivalent to Langlands functoriality.

His formulation and approach shows how a theory of $L$--functions can play a direct role in suggesting what the sought after [BK,BK2,FLN,Ge] Fourier transform and Poisson summation formula should be.
Moreover, the introduction of functions of $L$--type agrees completely with the approach of Braverman--Kazhdan [BK] and coincides with those of Godement--Jacquet [GJ] when $\rho$ is the standard representation of $GL_n(\bC)$.
It is therefore very tantalizing to see how different established approaches to the theory of automorphic $L$--functions themselves can contribute to our understanding of the Fourier transform and its Poisson summation formula in question.

Other approaches to the problem of the existence of the Fourier transform and its Poisson summation formula have been suggested, among them that of intertwining operators and Eisenstein series within the doubling method of Piatetski--Shapiro and Rallis which generalize that of Godement--Jacquet to classical groups, appearing in the recent manuscript of Wen--Wei Li [Li], and the preprint of Jacquet [J2] on symmetric square $L$--functions for $GL(2)$.
It remains to be seen how these different approaches will finally provide us with an answer to our question, one which is equivalent to functoriality and thus Langlands approach ``Beyond Endoscopy'' which also requires a Fourier transform and Poisson summation formula.

\parskip=0pt

\bigskip\bigskip\noindent
Department of Mathematics

\noindent
Purdue University

\noindent
150 N.~University Street

\noindent
West Lafayette, IN\ \ 47907

\noindent
email:\ \ shahidi\@math.purdue.edu

\bye